\documentclass[11pt,a4paper]{article}
\usepackage{amsfonts}
\usepackage{amsmath,amssymb}
\usepackage{mathrsfs}
\usepackage{hyperref}
\usepackage{graphicx}
\usepackage{pdfsync}

\newtheorem{thm}{Theorem}[section]
\newtheorem{cor}[thm]{Corollary}
\newtheorem{lem}[thm]{Lemma}
\newtheorem{prop}[thm]{Proposition}

\newtheorem{rem}[thm]{Remark}

\numberwithin{equation}{section}\allowdisplaybreaks
\newcommand{\N}{\ensuremath{\mathbb{N}}}
\newcommand{\R}{\ensuremath{\mathbb{R}}}

\newcommand{\X}{\mathfrak{X}_T^p}

\newcommand{\etu}{{T (u_0)}}
\newcommand{\CC}{{M^{\sigma_p}_c}}
\newcommand{\bph}{{\dot{B}^{-1+6/p}_{p/2, \infty}}}
\newcommand{\etl}{e^{t\Delta} }
\newcommand{\rt}{{\R^3}}
\newcommand{\lp}{{L^p}}
\newcommand{\e}{\epsilon }
\newcommand{\lineapp}{{e^{t\Delta}\big(\sum_{j\in \mathbb{J}_1^c(J)}\Lambda_{j,n}^{3/p}\phi_j + \psi_n^J \big)}}
\newcommand{\lmds}{\Lambda_{j,n}^{d/p}}
\newcommand{\lmdst}{\Lambda_{j,n}^{3/p}}
\newcommand{\lmdp}[1]{\Lambda_{#1,n}^{3/p}}
\newcommand{\sjjz}[1]{\sum_{j \in \mathbb{J}_{#1}(J)}}
\newcommand{\sjjzc}{\sum_{j \in \mathbb{J}_1^c(J)}}
\newcommand{\epn}{\epsilon_n}
\newcommand{\sgp}{\sigma_p}
\newcommand{\mtn}{T (u_{0,n})}
\newcommand{\U}{U_n^{app,J}}
\newcommand{\wt}{\widetilde T}
\newcommand{\re}{R_n^J}
\newcommand{\sol}[1]{NS(\Lambda^{3/p}_{#1,n} \phi_{#1})}
\newcommand{\al}{5p/6}
\newcommand{\be}{\beta(p)}

\newcommand{\fp}{3/p}

\newcommand{\smp}{\sigma(p)}
\newcommand{\ha}{1/2}
\newcommand{\ltb}{L_T^{\beta(p)}}

\newcommand{\anj}{A^J_n}
\newcommand{\bnj}{B^J_n}
\newcommand{\rnj}{R^J_n}

\def\leq{\leqslant}
\def\ge{\geqslant}

\def\leq{\leqslant}
\def\geq{\geqslant}

\begin{document}

\title{\Large\bf   Dynamical Behavior for the Solutions\\ of the Navier-Stokes Equation  }

\author{\normalsize \bf Kuijie Li$^{\dag}$, \ \ Tohru Ozawa$^{\ddag}$, \ \   Baoxiang Wang$^{\dag},$\footnote{Corresponding Author}   \\
\footnotesize
\it $^\dag$LMAM, School of Mathematical Sciences, Peking University, Beijing 100871, China, \\
\footnotesize
\it $^\ddag$ Department of Applied Physics, Waseda University, Tokyo 169-8555, Japan \\
\footnotesize
\it {Emails: kuijiel@pku.edu.cn, \ txozawa@waseda.jp, \ wbx@math.pku.edu.cn }\\
 } \maketitle

\thispagestyle{empty}
%\begin{multicols}{2}
\begin{abstract}

We study the Cauchy problem for the incompressible
Navier-Stokes equations (NS) in three and higher spatial dimensions:
\begin{align}
u_t -\Delta u+u\cdot \nabla u +\nabla p=0, \ \  {\rm div} u=0, \ \  u(0,x)=   u_0(x). \label{NSa}
\end{align}
Leray \cite{Le34} and Giga \cite{Gi86} obtained that for the weak and mild solutions $u$ of NS in $L^p(\mathbb{R}^d)$ which blow up at finite time $T>0$, respectively, one has that for $d <p \leq \infty$,
$$
\|u(t)\|_p \gtrsim ( T-t )^{-(1-d/p)/2}, \ \ 0< t<T.
$$
We will obtain the blowup profile and the concentration phenomena in $L^p(\mathbb{R}^d)$ with $d\leq p\leq \infty$ for the blowup mild solution.  On the other hand,  if the Fourier support has the form  ${\rm supp} \ \widehat{u_0} \subset \{\xi\in \mathbb{R}^n: \xi_1\geq L \}$ and $
  \|u_0\|_{\infty}  \ll L$ for some $L >0$,
 then \eqref{NSa} has a unique global solution $u\in C(\mathbb{R}_+, L^\infty)$. Finally, if the blowup rate is of type I:
 $$
\|u(t)\|_p \sim ( T-t )^{-(1-d/p)/2}, \ for \ 0< t<T<\infty, \ d<p<\infty
$$
in 3 dimensional case, then we can obtain a minimal blowup solution $\Phi$ for which
$$
\inf \{\limsup_{t \to T}(T-t)^{(1-3/p)/2}\|u(t)\|_{L^p_x}: \ u\in C([0,T); L^p) \mbox{\ solves \eqref{NSa}}\}
$$
is attainable at some $\Phi \in L^\infty (0,T; \ \dot B^{-1+6/p}_{p/2,\infty})$. \\



\end{abstract}

\section{Introduction} \label{sect1}

\medskip

We study   the Cauchy problem of the  incompressible
Navier-Stokes equations (NS) with initial data in $L^\infty (\mathbb{R}^d)$, $d\ge 2$:
\begin{align}
u_t -\Delta u+u\cdot \nabla u +\nabla p=0, \ \  {\rm div} u=0, \ \  u(0,x)=   u_0(x),  \label{NS}
\end{align}
where   $u=(u_1,...,u_d)$ denotes the flow
velocity vector and $p(t,x)$ describes the scalar pressure.   $\nabla=(\partial_1,...,\partial_d)$, $\Delta=\partial^2_1+...+\partial^2_d$,
$u_0(x)=(u^0_1,...,u^0_d)$ is a given velocity with $\mathrm{div} \
u_0=0$.
It is easy to see that \eqref{NS} can be rewritten as the following
equivalent form:
\begin{align}
  u_t -\Delta u+\mathbb{P}\ \textrm{div}(u\otimes u)=0, \ \
u(0, x)=   u_0(x),
 \label{NSE}
\end{align}
where $\mathbb{P}=I-\nabla \Delta^{-1}{\rm div}$ is the
projection operator onto the divergence free vector fields. The solution $u$ of NS formally satisfies the conservation
\begin{align}
 \frac{1}{2} \|u(t)\|^2_{L^2(\mathbb{R}^d)} +\int^t_0 \|\nabla u(\tau)\|^2_{L^2(\mathbb{R}^d)} d\tau = \frac{1}{2} \|u_0\|^2_{L^2(\mathbb{R}^d)}.
 \label{conservation}
\end{align}
It is known that \eqref{NS} is essentially equivalent to the following integral equation:
\begin{align} \label{NSI}
u(t) =   e^{ t \Delta } u_0 +  \int^t_0  e^{ (t-\tau) \Delta } \mathbb{P}\ \textrm{div}(u\otimes u)(\tau) d\tau
\end{align}
and the solution of \eqref{NSI} is said to be a {\it mild solution.} Note that \eqref{NS} is scaling invariant in the following sense:
if $u$ solves \eqref{NS},  so does   $u_{\lambda}(t, x)=\lambda u(\lambda^{2 }t, \lambda x)$ and $p_\lambda (t, x)=\lambda^2 p(\lambda^{2 }t, \lambda x)$ with initial data $\lambda u_0(\lambda x)$.
 A function space $X$ defined in $\mathbb{R}^d$ is said to be a \textit{ critical space} for \eqref{NS}
if the  norms of $u_{\lambda}(0,x)$ in $X$ are equivalent for
all $\lambda>0$ (i.e., $\|u_\lambda(0,\cdot)\|_X\sim \|u_0\|_X$). It is easy to see that $L^d$ and $\dot H^{d/2-1}$ are critical spaces of NS.

For the sake of convenience, we will denote by $NS(u_0)$ the solution of \eqref{NS} (or simply denote it by $u$ if there is no confusion), and by $T(u_0)$ the supremum of all $T>0$ so that the solution $NS(u_0)$ exists in the time interval $[0, T]$.

 Many years ago,  Leray \cite{Le34} showed that NS in 3D has at least one weak solution and he mentioned certain necessary blowup conditions for the weak solutions:
\begin{align} \label{blowup-Lp}
 \|NS(u_0)(t)\|_{L^p }  \gtrsim (T(u_0) -t)^{-( 1-d/p)/2 }, \ d < p \leq \infty.
\end{align}
The existence of the mild solution in $L^p$ was established by Kato in \cite{ Ka84} and the blowup rate \eqref{blowup-Lp} in all spatial dimensions  was recovered by Giga \cite{Gi86} for the mild solution in $C([0,T(u_0)); L^p)$ with $d<p<\infty$. The blowup in the critical space $L^3(\mathbb{R}^3)$  was first considered by  Escauriaza,  Seregin and  Sverak \cite{EsSeSv03} and they proved $\limsup_{t\to T(u_0)} \|NS(u_0)(t)\|_{L^3} $ $=\infty$ if $T(u_0)<\infty$, similar results in critical spaces $\dot H^{1/2}$ and $L^3$ were obtained in \cite{KeKo11, GaKoPl13} via the profile decomposition arguments developed by Kenig and Merle \cite{KeMe08} together with the backward uniqueness in \cite{EsSeSv03}. Seregin \cite{Se12} proved that $\lim_{t\to T(u_0)} \|NS(u_0)(t)\|_{L^3} =\infty$ if $T(u_0)<\infty$. Blowup results in $L^d$ for higher spatial dimensions were obtained in Dong and Du \cite{DoDu09} by following the approach in \cite{EsSeSv03}.
Recently, some generalizations for the blowup rates in $\dot H^s$ $(3/2\leq s\leq 5/2)$ were obtained in \cite{RoSaSi12, CoMoPi14},
\begin{align} \label{blowup-Hs}
 \|NS(u_0)(t)\|_{\dot H^s} \geq  c  (T(u_0)-t)^{-(  s  -1/2 )/2  }, \ \ d=3.
\end{align}
Some other kind of blowing up criteria can be found in Kozono, Ogawa and Taniuchi \cite{KoOg02} and references therein.

On the other hand, there are some works which have been devoted to generalize the initial data in some larger spaces;  cf. \cite{BaBiTa12, Can97,  Ch99, ChGaPa11, GiMi85, GiMi89, Iw10, KoTa01, Pl96, We80} and references therein. For the initial data in critical Besov type spaces, Cannone \cite{Can97}, Planchon \cite{Pl96}  and  Chemin \cite{Ch99} obtained global solutions in 3D for small data  in critical Besov spaces $\dot{B}^{3/p-1}_{p, q}$ for all $p<\infty, \ q\leq \infty$.
Bourgain and Pavlovic \cite{BoPa08} showed the  \textit{ill-posedness} (i.e., the solution map $u_0\to u$ is discontinuous) of NS in $\dot{B}^{-1}_{\infty,\infty}$,   Germain \cite{Ge08} proved that the solution map of NS is not $C^2$ in $\dot{B}^{-1}_{\infty, q} $ for any $q>2$,  Yoneda \cite{Yo10} showed that the solution map is discontinuous in $\dot{B}^{-1}_{\infty, q} $ for any $q>2$, and Wang \cite{Wa15} finally proved that NS is ill-posed in critical Besov spaces $\dot{B}^{-1}_{\infty, q} $, $1\leq q\leq 2$. So, NS is ill-posed in all critical Besov spaces $\dot{B}^{-1}_{\infty, q} $, $1\leq q\leq \infty$. Up to now, noticing the embedding  $\dot B^{-1+d/p}_{p,q} \subset BMO^{-1}$ $(p<\infty)$,  the known largest critical space for which NS is globally well posed for small initial data is $BMO^{-1}$, see Koch and Tataru \cite{KoTa01}.

Before stating our main result, we first give some notations. $C\geq 1, \ c\leq 1$ will denote constants which can be different at different places, we will use $A\lesssim B$ to denote   $A\leqslant CB$,  $A\sim B$ means that $A\lesssim B$ and $B\lesssim A$. We denote    by $L^p=L^p(\mathbb{R}^d)$ the Lebesgue space on which the norm is written as $\|\cdot\|_p$. $\|f\|_{\dot H^s} = \|(-\Delta)^{s/2}f\|_2$ and $H^s=L^2\cap \dot H^s$ for $s\geq 0$.
Let us write for any $\rho>0$,
\begin{align} \label{Rrho}
\mathbb{R}^d_{+,\rho} := \{\xi=(\xi_1,...,\xi_d) \in \mathbb{R}^d: \ \xi_1 \ge \rho\}.
\end{align}
The standard iteration sequence for NS is defined in the following way:
\begin{align} \label{iterationsequence}
u^{(n+1)}(t) =   e^{ t \Delta } u_0 -  \int^t_0  e^{ (t-\tau) \Delta } \mathbb{P}\ \textrm{div}(u^{(n)}\otimes u^{(n)})(\tau) d\tau, \ \ u^{(0)}=0.
 \end{align}
We will mainly consider the concentration behavior of the blowup solutions and obtain a global well-posedness result in $L^\infty$. The well-posedness of NS in $L^p$ with $d<p<\infty$ was established in Giga \cite{Gi86}:
\begin{thm}\label{th:1}
Let $u_0 \in L^p(\R^d)$ ($d < p <\infty$) be a divergence free vector.  Then there exists a time $T$ and a unique solution $NS(u_0)$ such that
the solution belongs to $ \mathfrak{X}^p_{T(u_0)}: = C([0,T(u_0)),L^p_x) \cap L^{p(1+2/d)}_{\rm loc}((0,T(u_0)),L^{p(1+2/d)}_x)$ and there exists $c_0>0$ such that $  u(t) $ blows up at $T(u_0)$ and for any $0\leq t< T(u_0)$,
\begin{equation}\label{lowerbound}
(T(u_0) -t)\| u(t)\|^{\sigma_p}_{p} \geq c_0,\, \,    \sigma _p := \frac {2}{  1-d/p }.
\end{equation}
\end{thm}
In Theorem \ref{th:1} the left case is $p=\infty$.  The Cauchy problem of the Navier-Stokes equations in $L^\infty$ and in $BUC$ spaces is studied by Cannone, Meyer \cite{Ca95,CaMe95}, Giga et al.
\cite{GiInMa99} and they proved a unique existence of a local-in-time solution in $L^\infty$ and in $BUC$ spaces. In \cite{GiInMa99}, the authors also obtained the smoothness of the solutions.     We will obtain a concentration phenomena of the blowup solutions with initial data only in $L^\infty$.  The first main result of this paper is
\begin{thm}\label{NS blowup and existence}
Let   $u_0 \in  L^\infty$, ${\rm div} u_0=0$. Then there exists a $T(u_0) >0 $ such that \eqref{NSI} has a unique solution $NS(u_0) \in L^\infty([0, T(u_0))\times \mathbb{R}^d)$. If  $T(u_0)  <\infty $, then
\begin{align} \label{blowuprate}
\omega(t):= \|NS(u_0) (t)\|_{ \infty } \gtrsim  (T(u_0) -t)^{- 1/2}, \ 0 \leq t< T(u_0).
\end{align}
For any $1\leq p \leq\infty$,  there exist  $x_n\in \mathbb{R}^d$ and $t_n \nearrow T(u_0) <\infty$  such that
\begin{align} \label{concentration}
 \|NS(u_0) (t_n)\|_{L^p(|x_n-\cdot| \lesssim \, \omega(t_n)^{-1})} \gtrsim  \omega(t_n)^{ 1-d/p }.
\end{align}
Moreover,  if ${\rm supp} \ \widehat{u_0} \subset \mathbb{R}^d_{+,\rho}$ for some $\rho>0$ and there exists $n_0 \in \mathbb{N}\cup \{0\}$ such that
\begin{align} \label{n0}
 4C(\|u^{(n_0)}\|_{L^\infty_{x, t\in [0,\infty)}} + \|u^{(n_0+1)}\|_{L^\infty_{x, t\in [0,\infty)}} ) \leq (n_0+1) \rho
 \end{align}
for some $C\gg 1$,  then \eqref{NSI} has a unique global solution $NS(u_0) \in C(\mathbb{R}_+, L^\infty)$.
\end{thm}

Let $\{x_{j,n}\} \subset \mathbb{R}^d$ and $\{\lambda_{j,n}\} \subset (0,\infty)$ be two sequences. $(\lambda_{j,n}, x_{j,n})^\infty_{n=1}$ $(j\in \mathbb{N})$ are said to be orthogonal sequences of scales and cores, if for any $j_1\neq j_2$, $j_1, j_2\in \mathbb{N}$, one has that
\begin{equation}\label{orthoa}
\lim_{n \to \infty} \bigg( \frac{\lambda_{j,n}}{\lambda_{k,n}} + \frac{\lambda_{k,n}}{\lambda_{j,n}}+
\frac{|x_{j,n} - x_{k,n}| }{\lambda_{j,n}} \bigg ) = \infty.
\end {equation}
It is known that $\lim \sup_{t\to T(u_0)} \|NS(u_0) (t)\|_{\dot H^{d/2-1}}=\infty$ if  $T(u_0)<\infty$ (cf. \cite{EsSeSv03, KeKo11, GaKoPl13, Se12, DoDu09}). In the next result we describe the blowup profiles for the blowing up solutions.

\begin{thm}\label{Blowup profile}

Let $u_0 \in H^{d/2-1}$, ${\rm div} u_0=0$ and $NS(u_0)\in C([0,T(u_0)); H^{d/2-1})$ be the solution of \eqref{NSI} which blows up at $T(u_0) <\infty$. Let  $t_n \nearrow T(u_0)$ satisfy
$$
\lim_{n\to \infty}\|NS(u_0) (t_n)\|_{\dot H^{d/2-1}}=\infty.
$$
Then there exist $\{\alpha_n\}$, $\{\phi_j\}$, $\{\lambda_j(t_n)\}$ and $\{x_{j}(t_n)\}$ with   $ \lim_{n\to \infty}\alpha_n = +\infty$, $\phi_j \in \dot H^{d/2-1}$, $\lim_{n\to \infty} \lambda_j(t_n) = 0$, $x_{j}(t_n)\in \mathbb{R}^d$, $(\lambda_j(t_n)\}, \{x_{j}(t_n))$ are orthogonal, such that  $NS(u_0) (t_n)$ can be decomposed into the following profiles:
\begin{align} \label{blowupdecompostion1}
NS(u_0) (t_n)= \alpha_n \left(\sum^J_{j=1} \frac{1}{\lambda_{j}(t_n)} \phi_j \left( \frac{x- x_{j}(t_n)}{\lambda_{j}(t_n)}  \right) + r^J_n \right),
\end{align}
where $r^J_n$ is a reminder that satisfies $\lim_{J\to \infty}\lim\sup_{n\to\infty}\|r^J_n\|_d=0$,  moreover,  we have $\alpha_n = \|NS(u_0) (t_n)\|_{\dot H^{d/2-1}}$, $\lambda_j (t_n) \lesssim  \alpha_n^{-2/(d-2)} $ and in particular, for $p\geq 2,$
\begin{align} \label{concentrationpicture1}
\sup_{\xi \in \mathbb{R}^d}  \|NS(u_0) (t_n)\|_{L^p(|\cdot- \, \xi| \lesssim   \lambda_{j}(t_n))}  \gtrsim   \lambda_{j}(t_n)^{d/p-1}   \|NS(u_0) (t_n)\|_{\dot H^{d/2-1}}.
\end{align}
\end{thm}

\begin{rem} \rm Theorems \ref{NS blowup and existence} and \ref{Blowup profile} needs several remarks.
\begin{itemize}
\item[\rm (i)]
 Noticing that $\omega (t) \geq c \  (T(u_0)-t)^{-1/2}$, we have for any $p\ge d$,
\begin{align} \label{concentration}
 \|NS(u_0) (t_n)\|_{L^p(|x_n-\cdot| \lesssim (T(u_0) -t_n)^{1/2}  )} \gtrsim   (T(u_0) -t_n)^{- (1  -d/ p)/2 }.
\end{align}

Taking $p=d$ in \eqref{concentration}, we find that
\begin{align} \label{concentration p=d}
 \|NS(u_0) (t_n)\|_{L^d(|x_n-\cdot| \lesssim (T(u_0) -t_n)^{1/2}  )} \gtrsim   1,
\end{align}
which implies that the solution has a concentration in a very small ball with radius less than or equals to $C\sqrt{T(u_0)-t_n}$ in $L^d$.

\item[\rm (ii)]   Taking $p=d$ in \eqref{concentrationpicture1} and noticing that $\lim_{n\to \infty}\|NS(u_0) (t_n)\|_{\dot H^{d/2-1}}=\infty$, we have
\begin{align} \label{concentrationpicture2}
\sup_{\xi \in \mathbb{R}^d}  \|NS(u_0) (t_n)\|_{L^d(|\cdot- \, \xi| \lesssim   \lambda_{j}(t_n))}  \gtrsim     \|u(t_n)\|_{\dot H^{d/2-1}},
\end{align}
which means that a very large potential norm is concentrated in a very small ball with radius less than or equals to $ C\|NS(u_0) (t_n)\|^{-2/(d-2)}_{\dot H^{d/2-1}}$. However, up to now, it is not very clear for us how to unify the concentration phenomena of \eqref{concentration p=d} and \eqref{concentrationpicture2}.

\item[\rm (iii)] In the blowup profile decomposition \eqref{blowupdecompostion1}, noticing that $\lambda_j(t_n)\to 0$ as $t_n \nearrow T(u_0)$,  we see that concentration blowup is the only way in all of the possible blowing up manners.

\item[\rm (iv)] Taking $n_0=0$ in \eqref{n0}, we see that condition \eqref{n0} can be substituted by the following condition:
\begin{align} \label{n00}
 4C \|u_0\|_{\infty}   \leq   \rho.
 \end{align}
Noticing that $\widehat{u_0}$ is supported in  $\mathbb{R}_{+,\rho}$, we see that condition \eqref{n00} contains a class of large data in $L^\infty$ if $\rho\gg 1$ which are out of the control of $\dot B^{-1+d/p}_{p,q}$ with $p<\infty$.

\end{itemize}
\end{rem}

\noindent Following \cite{EsSeSv03, KeKo11, GaKoPl13}, we see that, for the initial data $u_0$ in the critical spaces $X= \dot H^{1/2}(\mathbb{R}^3), \ L^3(\mathbb{R}^3)$, the solution $NS(u_0)$ enjoys the property
\begin{equation}\label{nonblowup}
\sup_{0<t< T (u_0)} \|NS(u_0)\|_{X} < \infty \ \ \Longrightarrow \ \  T(u_0) = \infty.
\end{equation}
We can further ask what happens if $X$ is a sub-critical spaces, say $X= \dot H^s(\mathbb{R}^3), \ s>1/2$, or $X=L^p(\mathbb{R}^3), \ p>3$. Noticing the blowup rate as in \eqref{blowup-Hs}, we have the following question in $H^s$:
$$
T(u_0) < \infty\,\, \Longrightarrow \limsup_{t \to T(u_0)} (T(u_0) -t)\|NS(u_0)(t)\|^{\sigma_s}_{\dot H^s} = \infty ?  \   (\sigma_s = \frac{2}{s-1/2})
$$
 E.~Poulon \cite{Po15} first considered such a kind of question and she introduced
$$
M^{\sigma_s}_{c}:= \inf_{u_0 \in \dot H^s \atop T(u_0) < \infty}\{\limsup_{t \to T(u_0)}(T(u_0)-t)\|NS(u_0)(t)\|^{\sigma_s}_{\dot H^s}\}.
$$
If $M_c^{\sigma_s}<\infty$, E.~Poulon proved that $M_c^{\sigma_s}$ can be attainable for some $u_0\in \dot H^s$ and the corresponding solution is uniformly bounded in critical Besov space $\dot B^{1/2}_{2,\infty}(\rt)$.  Following Theorem \ref{th:1}, we can define similar critical ``minimal'' quantity  adapted to  the $L^p(\rt)$ scale for $3<p<\infty$ ($\sigma_p =2/(1-3/p)$):
$$
M^{\sigma_p}_{c}:=\inf\left\{\limsup_{t \to T(u_0)}(T(u_0)-t)\|NS(u_0)\|^{\sigma_p}_{p} : u_0 \in {L^p}, \  T(u_0) < \infty \right\}.
$$
If $M^{\sigma_p}_{c}<\infty$, Giga's Theorem \ref{th:1} implies that there exists $u_0\in L^p$ such that the solution of NS blows up at finite time  $T(u_0)$ with the blowing up rate
$$
 \|NS(u_0)(t)\|_{p }   \sim (T(u_0) -t)^{-( 1-3/p)/2 }, \ 0\leq t< T(u_0).
$$
Such a kind of solution is said to be type-I blowing up solution (cf. \cite{KoNaSeSv09}). We have the following result:
\begin{thm} \label{mainthm}
Let $d=3, \ 3<p<\infty$, $\CC <\infty$.
Then there exists a ~$\Phi_0 \in L^p \cap \dot B^{-1+6/p}_{p/2,\infty}$ such that $\Phi:=NS(\Phi_0)$ blows up at time~$1$, and satisfies \begin{equation}\label{supcrit}
\sup_{0<\tau < 1 }(1-\tau)\|\Phi(\tau)\|^{\sigma_p}_{p} = \limsup_{\tau \to 1}(1-\tau)\|\Phi(\tau)\|^{\sigma_p}_{p} = {\CC}.
\end{equation}
Moreover, $\Phi$ lies in Besov spaces $\dot B^{-1+6/p}_{p/2,\infty}$, and
\begin{equation}\label{besov}
\sup_{0<\tau<1}\|\Phi(\tau)\|_{\dot B^{-1+6/p}_{p/2,\infty}} < \infty .
\end{equation}
\end{thm}

\subsection{Besov spaces}
 Let $\psi: \mathbb{R}^d\rightarrow
[0, 1]$ be a smooth cut-off function which equals $1$ on the closed ball $B(0,5/4):=\{\xi\in \mathbb{R}^d: |\xi|\le 5/4\}$ and equals $0$ outside the ball $B(0,3/2)$. Write
\begin{align} \label{phi}
\varphi(\xi):=\psi(\xi)-\psi(2\xi),  \ \ \varphi_j(\xi)=\varphi(2^{-j}\xi),
\end{align}
$\triangle_j:=\mathscr{F}^{-1}\varphi_j \mathscr{F}, \  j\in \mathbb{Z}$
are said to be the dyadic decomposition operators. One easily sees that ${\rm supp} \varphi_j \subset B(0, 2^{j+1})\setminus B(0, 2^{j-1})$. For convenience, we denote
\begin{align} \label{phij}
 P_{\le M}f := \mathscr{F}^{-1} \psi( \cdot \ /M ) \mathscr{F},  \ \ \ P_{\ge M} := I- P_{\le M}f.
\end{align}
 The norms in homogeneous Besov  spaces are defined as follows:
\begin{align}
\label{Besov} \|f\|_{\dot{B}^s_{p, q}}=
 \left(\sum_{j=-\infty}^{+\infty}2^{jsq}\|\triangle_j f\|^q_{p}\right)^{1/q}, \ \   \|f\|_{\dot{B}^s_{p, \infty}} =  \sup_{j\in \mathbb{Z}} 2^{jsq}\|\triangle_j f\|_{p}.
\end{align}
Using the heat kernel, we have (see \cite{BaChDa11,Tr83})
\begin{equation}\label{heatbesov}
\|u\|_{\dot B^s_{p,q}} \sim  \big \|t^{-\frac{s}{2}}\|\etl u\|_{p}\big\|_{L^q(\R^{+},\frac{dt}{t})} ,\quad  s<0.
\end{equation}

The rest of this paper is organized as follows. In Section~\ref{datainLinfty}, we consider the well-posedness and blowup concentration of NS in $L^\infty$ and prove Theorem \ref{NS blowup and existence}.  Using the  ``profile decomposition'' techniques, in Section \ref{blowupprof} we consider the blowup profile for the blowing up solution in $H^{d/2-1}$  and show Theorem \ref{Blowup profile}. In Section~\ref{pfmainthm} we will prove our Theorem~\ref{mainthm}, whose proof consists of two steps, constructing a critical solution in $\lp$ and $\bph$ respectively. The proof of Theorem~\ref{mainthm} relies upon a profile structure theorem, whose proof will be given in Sections~\ref{pfstructhm} and~\ref{pfthfour}. Finally, in the Appendix, we state some basic estimates on NS and prove a perturbation result which is useful in obtaining the estimate of the remainder term in the profile structure theorem.

\section{Initial data in $L^\infty$} \label{datainLinfty}

\subsection{Local well-posedness and blowup analysis}

We will frequently use the following Bernstein's multiplier estimate (see \cite{BL, WaHuHaGu11}):
\begin{lem}\label{Bern} {\rm (Multiplier estimate)} Let $L\in \mathbb{N}, \  L>n/2, \ \theta=n/2L$. We have
\begin{align}
\label{1.7}
\|\mathscr{F}^{-1} \rho\|_{1} \lesssim \|\rho\|^{1-\theta}_{2} \left(\sum^d_{i=1} \|\partial^L_i \rho\|_2^{\theta} \right).
\end{align}
\end{lem}
Recall that (see \cite{Ch99})
\begin{align} \label{exp-decay 0}
\left\|\triangle_j   e^{t\Delta} f \right\|_{r} \lesssim    e^{-t2^{2j-4}}
\| f\|_r, \ \ j\in \mathbb{Z}.
\end{align}
Similarly, we have
\begin{lem}\label{Exponential decay} {\rm (Exponential decay)} Let $1\le r\le \infty$.  We have
\begin{align}
\label{exp-decay}
\left\|\triangle_j (-\Delta)^{-1} \partial_\lambda \partial_\mu \partial_\nu  e^{t\Delta} f \right\|_{r} \lesssim   2^j e^{-t2^{2j-4}}
\| f\|_r, \ \ j\in \mathbb{Z}, \ \ 1\le\lambda, \mu, \nu \le d.
\end{align}
\end{lem}
{\bf Proof.} The idea follows from \cite{Ch99} (see also \cite{WaHuHaGu11}).  By Young's inequality, we have
\begin{align}
\label{exp-decay2}
\left\|\triangle_j (-\Delta)^{-1} \partial_\lambda \partial_\mu \partial_\nu  f \right\|_{r} \leq  \left\| \mathscr{F}^{-1} \left( \varphi_j \frac{\xi_{\lambda}\xi_\mu \xi_\nu}{|\xi|^2} e^{-t|\xi|^2} \right)  \right\|_{1}   \| f\|_r
\end{align}
By scaling argument and Lemma \ref{Bern}, we have
\begin{align}\label{exp-decay3}
\left\| \mathscr{F}^{-1} \left( \varphi_j \frac{\xi_{\lambda}\xi_\mu \xi_\nu}{|\xi|^2} e^{-t|\xi|^2} \right)  \right\|_{1} =  2^j \left\| \mathscr{F}^{-1} \left( \varphi  \frac{\xi_{\lambda}\xi_\mu \xi_\nu}{|\xi|^2} e^{-t2^{2j}|\xi|^2} \right)\right\|_1  \lesssim 2^j e^{-t2^{2j-4}}.
\end{align}
In view of \eqref{exp-decay2} and \eqref{exp-decay3}, we immediately have \eqref{exp-decay}.  $\hfill\Box$

For convenience, we denote
\begin{align}\label{inhomogeneous part}
  (\mathscr{A}_{t_0}  f)(t) := \int^t_{t_0} e^{(t-\tau)\Delta} f(\tau) d\tau.
\end{align}

\begin{lem}\label{Decay of higher frequency} {\rm (Decay of higher frequency)} Assume that  ${\rm supp} \ \widehat{f} \subset \{\xi : \ |\xi| \ge 2^{j_0}\}$.  Then we have
\begin{align}
\label{decay-higher-frequency}
\left\| \mathscr{A}_{t_0} \mathbb{P} \nabla f   \right\|_{L^\infty([t_0, \infty)\times \mathbb{R}^d)} \lesssim   2^{-j_0}
\| f\|_{L^\infty([t_0, \infty)\times \mathbb{R}^d)}.
\end{align}
\end{lem}
{\bf Proof.} By the dyadic decomposition,
\begin{align}
\label{decay-higher-frequency2}
  f = \sum_{j\ge j_0-1} \triangle_j f.
\end{align}
Using Lemma \ref{Exponential decay},  one sees that
\begin{align}\label{decay-higher-frequency3}
\left\| \mathscr{A}_{t_0} \mathbb{P} \nabla f   \right\|_{L^\infty([t_0, \infty)\times \mathbb{R}^d)}  & \leq   \sum_{j\ge j_0-1}  \left\| \triangle_j \mathscr{A}_{t_0} \mathbb{P} \nabla f   \right\|_{L^\infty([t_0, \infty)\times \mathbb{R}^d) } \nonumber\\
&  \lesssim   \sum_{j\ge j_0-1} 2^j {\rm sup_{t\ge t_0}} \int^t_{t_0}  e^{-(t-\tau)2^{2j-4}} \|f(\tau)\|_\infty d\tau \nonumber \\
&  \lesssim   \sum_{j\ge j_0-1} 2^{-j} {\rm sup_{t\ge t_0}} \left(1-  e^{-(t-t_0)2^{2j-4}} \right) \|f\|_{L^\infty([t_0, \infty)\times \mathbb{R}^d)} \nonumber\\
&  \lesssim     2^{-j_0}   \|f \|_{L^\infty([t_0, \infty)\times \mathbb{R}^d)}.
\end{align}
The result follows. $\hfill\Box$

\begin{lem}\label{Short time estimates of lower frequency} {\rm (Short time estimates of lower frequency)} Let $j_0\in \mathbb{Z}$.  We have for any $ f \in L^\infty ([t_0, t_1]\times \mathbb{R}^d)$,
\begin{align}
\label{lower-frequency-estimate}
\left\| P_{\le 2^{ j_0}}\mathscr{A}_{t_0} \mathbb{P} \nabla f   \right\|_{L^\infty([t_0, t_1]\times \mathbb{R}^d)} \lesssim   2^{j_0}(t_1- t_0)
\| f\|_{L^\infty([t_0, t_1]\times \mathbb{R}^d)}.
\end{align}
\end{lem}
{\bf Proof.} Since $e^{-t|\xi|^2} (\xi_{j} \xi_{k} \xi_{l} |\xi|^{-2}) \in L^1 (\mathbb{R}^d)$ for any $t>0$, whose Fourier transform has an integral form,  it follows that
\begin{align}
\label{without delta}
 \mathscr{F}^{-1} (e^{-t|\xi|^2}  \xi_{j} \xi_{k} \xi_{l} |\xi|^{-2})     =   \sum_{j\in \mathbb{Z}} \triangle_j   \mathscr{F}^{-1} (e^{-t|\xi|^2}  \xi_{j} \xi_{k} \xi_{l} |\xi|^{-2}).
\end{align}
By Lemma \ref{Exponential decay} and \eqref{without delta}, we have
\begin{align}
\label{lower-frequency-estimate2}
\left\| P_{\le 2^{j_0}}\mathscr{A}_{t_0} \mathbb{P} \nabla f   \right\|_{L^\infty([t_0, t_1]\times \mathbb{R}^d)}
 & \leq   \sum_{j\leq j_0+1}  \left\| \triangle_j \mathscr{A}_{t_0} \mathbb{P} \nabla f   \right\|_{L^\infty([t_0, t_1]\times \mathbb{R}^d) } \nonumber\\
 & \lesssim   \sum_{j\leq j_0+1} 2^j \int^{t_1}_{t_0}   \|f(\tau)\|_\infty d\tau  \nonumber\\
& \lesssim  2^{j_0}(t_1- t_0)
\| f\|_{L^\infty([t_0, t_1]\times \mathbb{R}^d)},
\end{align}
which is the result, as desired. $\hfill\Box$

\begin{lem}\label{Local well posedness} {\rm (Local well posedness)} Let $u_0\in L^\infty$ with ${\rm div} u_0=0$.  There exists a $T>0$ such that \eqref{NSI} has a unique solution $NS(u_0) \in L^\infty ([0,T)\times \mathbb{R}^d)$. If $T<\infty$, then we have  $\lim_{t\to T }\|NS(u_0) (t)\|_{ \infty } =\infty $ and
\eqref{blowuprate} holds true.
\end{lem}
{\bf Proof.} The proof of the local existence can be found in \cite{GiInMa99}. However, the blowup rate \eqref{blowuprate} is very important for our later purpose and we sketch the proof here.  Put
\begin{align}
\label{woking space}
\mathscr{D}= \{u\in L^\infty([0,t_0]\times \mathbb{R}^d):   \ \|u\|_{L^\infty([0,t_0]\times \mathbb{R}^d)} \le 2C \|u_0\|_{\infty} \}
\end{align}
with natural metric $d(u,v)= \|u-v\|_{L^\infty([0,t_0]\times \mathbb{R}^d)}$. Let us consider the mapping
\begin{align} \label{mapping}
\mathscr{T}: u(t) \to   e^{ t \Delta } u_0 +  \mathscr{A}_0 \mathbb{P}\ \textrm{div}(u\otimes u).
\end{align}
Using $e^{t\Delta}: L^\infty\to L^\infty$, and Lemmas \eqref{Decay of higher frequency} and \eqref{Short time estimates of lower frequency}, we have
\begin{align}
\label{mapping2}
\left\| \mathscr{T}u   \right\|_{L^\infty([0, t_0]\times \mathbb{R}^d)} \leq C \|u_0\|_\infty + C (2^{-j_0} + t_0 2^{j_0})  \|  u  \|^2_{L^\infty([0, t_0]\times \mathbb{R}^d)}.
\end{align}
Taking $j_0$ such that $2^{-2j_0} \sim t_0 $, one has that
\begin{align}
\label{mapping3}
\left\| \mathscr{T}u   \right\|_{L^\infty([0, t_0]\times \mathbb{R}^d)} \leq C \|u_0\|_\infty + C  \sqrt{t_0}  \|  u  \|^2_{L^\infty([0, t_0]\times \mathbb{R}^d)}.
\end{align}
Similarly,
\begin{align}
\label{mapping4}
\left\| \mathscr{T}u  - \mathscr{T}v  \right\|_{L^\infty([0, t_0]\times \mathbb{R}^d)} \leq   C  \sqrt{t_0} \left(  \sum_{w=u,v}\|w\|_{L^\infty([0, t_0]\times \mathbb{R}^d)} \right) \|u-v  \|_{L^\infty([0, t_0]\times \mathbb{R}^d)}.
\end{align}
Further, one can choose $t_0$ verifying
\begin{align}
\label{mapping5}
   \sqrt{t_0} = 1/(8 C^2\|u_0\|_\infty ).
\end{align}
We easily see that $\mathscr{T}$ is a contraction mapping from $\mathscr{D}$ into itself. So, $\mathscr{T}$ has a unique fixed point in $\mathscr{T}$, which is a solution of \eqref{NSI}. It is easy to see that $u\in C((0,T]; L^\infty)$ (see \cite{GiInMa99}, for instance).  The solution can be extended exactly in the same way as above. Indeed, considering the mapping:
\begin{align} \label{secondmapping}
\mathscr{T}_1: u(t) \to   e^{ (t-t_0) \Delta } u(t_0) +  \mathscr{A}_{t_0} \mathbb{P}\ \textrm{div}(u\otimes u),
\end{align}
one has that
\begin{align}
\label{secondmapping3}
\left\| \mathscr{T}_1 u   \right\|_{L^\infty([t_0, t_1]\times \mathbb{R}^d)} \leq C \|u (t_0)\|_\infty + C  \sqrt{t_1-t_0}  \|  u  \|^2_{L^\infty([t_0, t_1]\times \mathbb{R}^d)}.
\end{align}
Similarly,
\begin{align}
\label{secondmapping4}
& \left\| \mathscr{T}_1u  - \mathscr{T}_1v  \right\|_{L^\infty([t_0, t_1]\times \mathbb{R}^d)} \nonumber\\
& \ \ \ \ \leq   C  \sqrt{t_1-t_0} \left(  \sum_{w=u,v}\|w\|_{L^\infty([t_0, t_1]\times \mathbb{R}^d)} \right) \|u-v  \|_{L^\infty([t_0, t_1]\times \mathbb{R}^d)}.
\end{align}
So, we can extend the solution from $[0,t_0]$ to $[t_0, t_1]$ if
\begin{align}
\label{secondmapping5}
   \sqrt{t_1- t_0} = 1/(8 C^2\|u(t_0)\|_\infty ).
\end{align}
Repeating the procedure as above, we can extend the solution step by step to $[t_1,t_2], \ [t_2, t_3],...$ if
\begin{align}
\label{secondmapping6}
   \sqrt{t_i- t_{i-1}} = 1/(8 C^2\|u(t_{i-1})\|_\infty ), \ \ i=0,1,2,... \ .
\end{align}
Now let us assume that $t_i \nearrow T$.  If $T<\infty$, we easily see that $\limsup_{t\to T} \|NS(u_0)(t)\|_\infty= \infty$.  Moreover, we can show \eqref{blowuprate} holds true.  Assume for a contrary that there exist two sequences $s_n \nearrow T$ and $c_n \searrow 0 $ satisfying
\begin{align} \label{blowupratesequence}
  \|NS(u_0)(s_n)\|_{ \infty } \leq c_n (T -s_n)^{- 1/2}, \ \ n=1,2,... \ .
\end{align}
Observing the integral equation
\begin{align} \label{NSIblowup}
 u(t) =  e^{ (t-s_n) \Delta } NS(u_0) (s_n) +  \mathscr{A}_{s_n} \mathbb{P}\ \textrm{div}(u\otimes u),
\end{align}
similarly as in  \eqref{secondmapping3},  we have
 for any $s_n \leq s < T$,
\begin{align}
\label{NSIblowup2}
 \| u \|_{L^\infty([s_n, s]\times \mathbb{R}^d)}  & \leq C \|NS(u_0) (s_n)\|_\infty + C  \sqrt{s-s_n}  \|  u  \|^2_{L^\infty([s_n, s]\times \mathbb{R}^d)} \nonumber\\
    & \leq C c_n (s -s_n)^{- 1/2}  + C  \sqrt{s-s_n}  \|  u  \|^2_{L^\infty([s_n, s]\times \mathbb{R}^d)}.
\end{align}
\eqref{NSIblowup2} implies that
\begin{align}
\label{NSIblowup3}
  \sqrt{s -s_n} \| u\|_{L^\infty([s_n, s]\times \mathbb{R}^d)}
      \leq C c_n + C   (\sqrt{s-s_n}   \|  u  \|_{L^\infty([s_n, s]\times \mathbb{R}^d)})^2.
\end{align}
It follows from  \eqref{NSIblowup3} that
\begin{align}
\label{NSIblowup4}
  \sqrt{s -s_n} \| u\|_{L^\infty([s_n, s]\times \mathbb{R}^d)} \geq 1/2C,  \ \ or  \ \  \sqrt{s -s_n} \| u\|_{L^\infty([s_n, s]\times \mathbb{R}^d)} \leq  2Cc_n.
\end{align}
Taking $c_n \ll 1$ and $s\to T$,  we see that \eqref{NSIblowup4} contradicts with \eqref{blowupratesequence} or with the fact that $\limsup_{t\to T}\|NS(u_0) (t)\|_\infty =\infty$.  $\hfill\Box$

\begin{lem}\label{Concentration of local solutions} {\rm (Concentration)}    Let $T<\infty$  and $NS(u_0) \in L^\infty ([0,T)\times \mathbb{R}^d)$ be the solution of \eqref{NSI} obtained in Lemma \ref{Local well posedness}.  Then for any $1\le p<\infty$,  there exist two sequences $\{x_n\}$ and $\{\tau_n\}$ of $x_n\in \mathbb{R}^d$ and $\tau_n \nearrow T $ satisfying
\begin{align} \label{concentrationfinite}
 \|NS(u_0) (\tau_n)\|_{L^p(|x_n-\cdot| \lesssim \, \omega(\tau_n)^{-1})} \gtrsim  \omega(\tau_n)^{ 1-d/p }, \ \ \omega (t) = \|NS(u_0) (t)\|_\infty.
\end{align}
\end{lem}
{\bf Proof. } Put $\tau_0=0$. Since $\lim_{t\to T}\omega(t)=\infty$, we can find a $t_1>\tau_0$ such that
\begin{align} \label{concentrationproof1}
  \omega (t_1) \geq 200 C^2 \omega (\tau_0).
\end{align}
We can further find a $\tau_1\in (\tau_0, t_1]$ verifying
\begin{align} \label{concentrationproof2}
  \omega (\tau_1) \geq \frac{1}{2} \sup_{t\in [\tau_0, t_1]} \omega (t).
\end{align}
Then we have
\begin{align} \label{concentrationproof3}
  \omega (\tau_1) \geq \max\left\{ 100 C^2 \omega (\tau_0), \ \  \frac{1}{2} \sup_{t\in [\tau_0, \tau_1]} \omega (t) \right\}.
\end{align}
Repeating this procedure, one can choose a monotone sequence $\{\tau_n\}$ verifying
\begin{align} \label{concentrationproof3}
  \omega (\tau_n) \geq \max\left\{ 100 C^2 \omega (\tau_{n-1}), \ \  \frac{1}{2} \sup_{t\in [\tau_{n-1}, \tau_n]} \omega (t) \right\}.
\end{align}
Since $\omega (\tau_n) \to \infty$, we have $\tau_n \nearrow T$.

{\it Claim. } For simplicity, we write $u:=NS(u_0)$.  We have
\begin{align} \label{concentrationinlow}
\|P_{\leq 100 C  \omega(\tau_n)} u(\tau_n)\|_{\infty} \geq \frac{1}{2} \omega (\tau_n) .
\end{align}
In fact, if \eqref{concentrationinlow} does not hold, then one has that
\begin{align} \label{concentrationinhigh}
\|P_{> 100 C  \omega(\tau_n)} u(\tau_n)\|_{\infty} \geq \frac{1}{2} \omega (\tau_n) .
\end{align}
Let us consider the integral equation
\begin{align} \label{nth-mapping}
 u(t) =   e^{ (t-\tau_{n-1}) \Delta } u(\tau_{n-1}) +  \mathscr{A}_{\tau_{n-1}} \mathbb{P}\ \textrm{div}(u\otimes u).
\end{align}
It follows from Lemma \ref{Decay of higher frequency} and \eqref{concentrationinhigh} that
\begin{align} \label{nth-mapping2}
 \omega (\tau_n) & \leq 2 \|P_{> 100 C  \omega(\tau_n)} u(\tau_n)\|_{\infty} \nonumber\\
  & \leq 2 C   \| P_{> 100  C \omega(\tau_n)} u(\tau_{n-1})\|_\infty +  2\|P_{> 100  C \omega(\tau_n)}\mathscr{A}_{\tau_{n-1}} \mathbb{P}\ \textrm{div}(u\otimes u)\|_{L^\infty([\tau_{n-1}, \tau_n]\times \mathbb{R}^d)} \nonumber\\
 & \leq 2 C^2  \| u(\tau_{n-1})\|_\infty +  \frac{1}{50   \omega(\tau_n)}\|u \|^2_{L^\infty([\tau_{n-1}, \tau_n]\times \mathbb{R}^d)} \nonumber\\
 & \leq \frac{1}{50}  \omega(\tau_n) +  \frac{4}{50}\omega(\tau_n) <  \frac{1}{2}  \omega(\tau_n).
\end{align}
A contradiction! So, we have \eqref{concentrationinlow}. For convenience, we write
\begin{align} \label{notionbeta}
 \beta(t): = 100 C \omega (t)  =  100 C \|u(t)\|_\infty.
\end{align}
Now we prove \eqref{concentrationfinite}. There exists $x_n \in \mathbb{R}^d$ such that
\begin{align} \label{concentrationinlowpf1}
\|P_{\leq  \beta(\tau_n)} u(\tau_n)\|_{\infty} \leq & 2 \beta(\tau_n)^d \left| \int_{\mathbb{R}^d} (\mathscr{F}^{-1}{\psi}) (\beta(\tau_n)(x_n -y)) u(\tau_n, y) dy \right| \nonumber\\
\leq  & 2 \beta(\tau_n)^d \left| \int_{|x_n-y|\leq M \beta(\tau_n)^{-1} } (\mathscr{F}^{-1}{\psi}) (\beta(\tau_n)(x_n -y)) u(\tau_n, y) dy \right| \nonumber\\
& + 2 \beta(\tau_n)^d \left| \int_{|x_n-y|\geq M \beta(\tau_n)^{-1} } (\mathscr{F}^{-1}{\psi}) (\beta(\tau_n)(x_n -y)) u(\tau_n, y) dy \right| \nonumber\\
:= & \ I+II.
\end{align}
By H\"older's inequality,
\begin{align} \label{concentrationinlowpf2}
II  \leq  2 \|\mathscr{F}^{-1}{\psi} \|_{L^1(|\cdot| \geq M )}  \|u(\tau_n)\|_\infty \leq 4 \|\mathscr{F}^{-1}{\psi} \|_{L^1(|\cdot| \geq M )}  \|P_{\leq \beta(\tau_n)} u(\tau_n)\|_\infty.
\end{align}
Taking a suitable $M\gg 1$, one easily sees that
\begin{align} \label{concentrationinlowpf3}
II  \leq \frac{1}{4}  \|P_{\leq \beta(\tau_n)} u(\tau_n)\|_\infty.
\end{align}
So, it follows from \eqref{concentrationinlowpf3} and H\"older's inequality that
\begin{align} \label{concentrationinlowpf4}
\|P_{\leq  \beta(\tau_n)} u(\tau_n)\|_{\infty}
\leq  & 4 \beta(\tau_n)^d \left| \int_{|x_n-y|\leq M \beta(\tau_n)^{-1} } (\mathscr{F}^{-1}{\psi}) (\beta(\tau_n)(x_n -y)) u(\tau_n, y) dy \right| \nonumber\\
\lesssim   &   \beta(\tau_n)^{d/p} \|u(\tau_n)\|_{L^p (|x_n-\cdot|\leq M \beta(\tau_n)^{-1}) }.
\end{align}
So, From \eqref{concentrationinlow} and \eqref{concentrationinlowpf4} we immediately have the result, as desired. $\hfill \Box$

\subsection{Global well-posedness }

For convenience, we denote
\begin{align} \label{bilinearop}
\mathscr{B}(u,v) =      \int^t_0  e^{ (t-\tau) \Delta } \mathbb{P}\ \textrm{div}(u \otimes v)(\tau) d\tau.
\end{align}
By Lemma \ref{Decay of higher frequency}, we have

\begin{lem}\label{Gocal solutions 1} {\rm (Decay with frequency superposition)} Let  ${\rm supp} \ \widehat{u_i} \subset \mathbb{R}^d_{+,\rho_i}$, $i=1,2$. Then we have  ${\rm supp} \ \widehat{\mathscr{B}(u_1, u_2)} \subset \mathbb{R}^d_{+,\rho_1+ \rho_2}$ and
\begin{align} \label{superpositiondecay}
 \|\mathscr{B}(u_1,u_2)\|_{L^\infty_{x,t\in [0,\infty)}}  \leq   C(\rho_1+\rho_2)^{-1} \prod^2_{i=1}\| u_i\|_{L^\infty_{x,t\in [0,\infty)}}.
\end{align}
\end{lem}
Let us recall the iteration sequence
\begin{align} \label{iterationsequence1}
 & u^{(0)}=0, \ \ u^{(1)}(t) =   e^{ t \Delta } u_0,  \nonumber \\
 & u^{(n+1)}(t) =   e^{ t \Delta } u_0 -   \mathscr{B}(u^{(n)}, u^{(n)}), \ n=1,2,...\ .
 \end{align}

\begin{lem}\label{Gocal solutions 2} {\rm (Frequency superposition of higher iteration terms)} Let  ${\rm supp} \ \widehat{u_0} \subset \mathbb{R}^d_{+,\rho}$. Then we have
\begin{align} \label{superposition of higher terms}
{\rm supp} \ \widehat{u^{(n)}} \subset \mathbb{R}^d_{+,\rho}, \  \   {\rm supp} \ (u^{(n+1)}-  u^{(n)})^{\widehat{  }} \ \subset   \mathbb{R}^d_{+, (n+1)\rho}.
\end{align}
\end{lem}
{\bf Proof.} By \eqref{iterationsequence1}, ${\rm supp} \ \widehat{u_0} \subset \mathbb{R}^d_{+,\rho}$ and induction we see that  ${\rm supp} \ \widehat{u^{(n)}} \subset \mathbb{R}^d_{+,\rho}$. Let us observe that
\begin{align} \label{induction1}
     (u^{(n+1)}-  u^{(n)})  = \mathscr{B} (  u^{(n-1)} - u^{(n)} ,  \  u^{(n)} ) + \mathscr{B} ( u^{(n-1)}, \   u^{(n-1)}- u^{(n)} )
\end{align}
and
\begin{align} \label{induction2}
     (u^{(2)}-  u^{(1)})  = - \mathscr{B} ( e^{t\Delta}u_0, \   e^{t\Delta}u_0).
\end{align}
By \eqref{induction2} and ${\rm supp} \ \widehat{u_0} \subset \mathbb{R}^d_{+,\rho}$, we see that  $     {\rm supp} \ (u^{(2)}-  u^{(1)})^{\widehat{    }} \ \subset   \mathbb{R}^d_{+, 2\rho}.$  By ${\rm supp} \ \widehat{u^{(n)}} \subset \mathbb{R}^d_{+,\rho}$ and induction, it follows from \eqref{induction1} that  \eqref{superposition of higher terms} holds true. $\hfill\Box$

\begin{lem}\label{Iteration estimate} {\rm (Compactness of Iteration)} Let  ${\rm supp} \ \widehat{u_0} \subset \mathbb{R}^d_{+,\rho}$. Assume that there exists $n_0\in \mathbb{N}$ such that
\begin{align} \label{iterationestimate1}
M_0:=  \| u^{(n_0)}\|_{L^\infty_{x,t\in [0,\infty)}} + \| u^{(n_0+1)}\|_{L^\infty_{x,t\in [0,\infty)}}  \leq  \frac{(n_0+1) \rho}{8C}.
\end{align}
Then $\{u^{(n)}\}$  is a Cauchy sequence in $L^\infty_{x,t\in [0,\infty)}$.
\end{lem}
{\bf Proof.} By \eqref{induction1}, we have for any $n>n_0 $,
\begin{align} \label{compact1}
     \|u^{(n+1)}-  u^{(n)}\|_{L^\infty_{x,t\in [0,\infty)}}  = & \|\mathscr{B} (  u^{(n)}-  u^{(n-1)},  \  u^{(n)} )\|_{L^\infty_{x,t\in [0,\infty)}} \nonumber\\
     &  + \|\mathscr{B} ( u^{(n-1)}, \  u^{(n)}-  u^{(n-1)} )\|_{L^\infty_{x,t\in [0,\infty)}}.
\end{align}
By Lemmas \ref{Gocal solutions 1} and  \ref{Gocal solutions 2}, it follows from \eqref{compact1} that
\begin{align} \label{compact2}
       \|\mathscr{B} (  u^{(n)}-  u^{(n-1)},  \  u^{(n)} )\|_{L^\infty_{x,t\in [0,\infty)}}  \leq   \frac{C}{(n+1)\rho}
         \|  u^{(n)}\|_{L^\infty_{x,t\in [0,\infty)}}
       \|  u^{(n)}-  u^{(n-1)}\|_{L^\infty_{x,t\in [0,\infty)}}.
\end{align}
Using a similar way, one can estimate the second term in \eqref{compact1}. So, in view of \eqref{compact1} and \eqref{compact2} we have
\begin{align} \label{compact3}
     \|u^{(n+1)}-  u^{(n)}\|_{L^\infty_{x,t\in [0,\infty)}}   \leq  &  \frac{C}{(n+1)\rho}
         (\|  u^{(n)}\|_{L^\infty_{x,t\in [0,\infty)}} + \|  u^{(n-1)}\|_{L^\infty_{x,t\in [0,\infty)}}) \nonumber\\
         & \times  \|  u^{(n)}-  u^{(n-1)}\|_{L^\infty_{x,t\in [0,\infty)}}.
\end{align}
Repeating the procedure as in \eqref{compact3}, one can obtain that for any $n>n_0$,
\begin{align} \label{compact4}
     \|u^{(n+1)}-  u^{(n)}\|_{L^\infty_{x,t\in [0,\infty)}}   \leq  &  \frac{C^{n-n_0} (n_0+1)!}{\rho^{n-n_0} (n+1)!}
     \prod^{n-1}_{k=n_0}    (\|  u^{(k+1)}\|_{L^\infty_{x,t\in [0,\infty)}} + \|  u^{(k)}\|_{L^\infty_{x,t\in [0,\infty)}}) \nonumber\\
         & \times  \|  u^{(n_0+1)}-  u^{(n_0)}\|_{L^\infty_{x,t\in [0,\infty)}}.
\end{align}
Now we show the following

{\it Claim.}  $\{u^{(n)}\}_{n>n_0+1}$ is bounded in $L^\infty_{x,t\in [0,\infty)}$ and $\|u^{(n)}\|_{L^\infty_{x,t\in [0,\infty)}} \leq 2 M_0$, $n>n_0+1$.

We prove the Claim. Taking $n=n_0+1$ in \eqref{compact3}, we have from \eqref{iterationestimate1} that
\begin{align} \label{compact5}
     \|u^{(n_0+2)}-  u^{(n_0+1)}\|_{L^\infty_{x,t\in [0,\infty)}}   \leq    \frac{C M_0}{(n_0+2)\rho} M_0 \leq \frac{M_0}{2}.
\end{align}
Hence, we have $\|u^{(n_0+2)}\|_{L^\infty_{x,t\in [0,\infty)}} \leq 2 M_0$.  Now let us assume the following induction assumption holds true:
\begin{align} \label{compact6}
     \|u^{(m)}\|_{L^\infty_{x,t\in [0,\infty)}}   \leq    2M_0, \ \ n_0+2 \leq m \leq n.
\end{align}
We show that \eqref{compact6} also holds for $m=n+1$.  Applying \eqref{compact4}, one has that for any $\ell= n_0+1,...,n-1$,
\begin{align} \label{compact7}
     \|u^{(\ell+1)}-  u^{(\ell)}\|_{L^\infty_{x,t\in [0,\infty)}}   \leq  &  \frac{C^{\ell-n_0} (n_0+1)!}{\rho^{\ell-n_0} (\ell+1)!}
     \prod^{\ell-1}_{k=n_0}    (\|  u^{(k+1)}\|_{L^\infty_{x,t\in [0,\infty)}} + \|  u^{(k)}\|_{L^\infty_{x,t\in [0,\infty)}}) \nonumber\\
         & \times  \|  u^{(n_0+1)}-  u^{(n_0)}\|_{L^\infty_{x,t\in [0,\infty)}} \nonumber\\
     \leq  &  \frac{C^{\ell-n_0} (4M_0)^{\ell-n_0} (n_0+1)!}{\rho^{\ell-n_0} (\ell+1)!}  M_0  \leq  \frac{M_0}{2^{\ell-n_0}}.
\end{align}
We have
\begin{align} \label{compact8}
     \|u^{(n+1)}-  u^{(n_0+1)}\|_{L^\infty_{x,t\in [0,\infty)}}   & \leq  \sum^{n}_{\ell=n_0+1} \|u^{(\ell+1)}-  u^{(\ell)}\|_{L^\infty_{x,t\in [0,\infty)}} \nonumber\\
      & \leq  \sum^{n}_{\ell=n_0+1}     \frac{M_0}{2^{\ell-n_0}} \leq M_0.
\end{align}
It follows that $\|u^{(n+1)}\|_{L^\infty_{x,t\in [0,\infty)}} \leq 2M_0$.

Finally, we show that $\{u^{(n)}\}$ is a Cauchy sequence. Again, in view of \eqref{compact7}, we have for any $n>m \gg  n_0$
\begin{align} \label{compact9}
     \|u^{(n+1)}-  u^{(m+1)}\|_{L^\infty_{x,t\in [0,\infty)}}   & \leq  \sum^{n}_{\ell=m+1} \|u^{(\ell+1)}-  u^{(\ell)}\|_{L^\infty_{x,t\in [0,\infty)}} \nonumber\\
      & \leq  \sum^{n}_{\ell=m+1}     \frac{M_0}{2^{\ell-n_0}} \leq \frac{M_0}{2^{m-n_0}} \to 0, \ \ m\to \infty.
\end{align}
Therefore, $\{u^{(n)}\}$ is a Cauchy sequence in $L^\infty_{x,t\in [0,\infty)}$. $\hfill \Box$

\subsection{Some remarks on the global well-posedness}

Recall that in \cite{Can97, Pl96, Ch99, BoPa08,Ge08,Yo10,Wa15} it was shown that:

\begin{thm}\label{NSill}
Let $1\leq p< \infty, \ 1\leq q\leq \infty$. Then NS is globally well-posed in $\dot B^{-1+d/p}_{p,q}$ for sufficiently small data.  Moreover, NS is  ill posed in $\dot B^{-1}_{\infty,q}$, i.e., the solution map $u_0\to u$ is discontinuous.
\end{thm}

However, Theorem \ref{NS blowup and existence} implies the well-posedness for a class of large data in $\dot B^{-1+d/p}_{p, q}$.  Let $\varphi$ be as in \eqref{phi},  $\widetilde{\varphi}(\xi) = \varphi (\xi) \chi_{\{1/2\leq \cdot\leq 2 \}}(\xi)$\footnote{We denote by $\chi_E$ the characteristic function on $E$.} and
\begin{align}
\label{delta}  \widehat{u^{1}_0} (\xi) = c \, 2^{L(1-d)} \prod^d_{i=1}  \widetilde{\varphi} (2^{-L} \xi_i),     \ \  \widehat{u^{2}_0} =- \frac{\xi_1}{\xi_2} \widehat{u^{1}_0} (\xi), \ \ u_0^3=...=u_0^n=0,
\end{align}
where $c \ll 1$ is a small constant. It is easy to see that ${\rm div} u_0=0$. We have,
\begin{align}
\label{deltaexample}  \|u^{1}_0\|_\infty  \sim  2^{L}, \ \ \|u^{2}_0\|_\infty  \sim   2^{L}.
\end{align}
It follows from Theorem \ref{NS blowup and existence} that for the initial data as in \eqref{delta}, NS is globally well-posed in $L^\infty$.  Also, we see that $\|u_0\|_{\dot B^{-1}_{\infty, q}} \ll 1$ and $\|u_0\|_{\dot B^{-1+d/p}_{p, q}} \sim 2^{dL/p}$, which is arbitrarily large in $\dot B^{-1+d/p}_{p, q}$.\\

\section{Blowup profile} \label{blowupprof}

If we consider the profile decomposition in $\dot H^s$, the orthogonal condition \eqref{orthoa} can be replaced by
\begin{align} \label{orthogonal}
{\it either } \ \ \lim_{n\to \infty} \left|\ln\frac{\lambda_{j_1,n}}{\lambda_{j_2,n}} \right|   = \infty, \ \ {\it or } \ \  \frac{\lambda_{j_1,n}}{\lambda_{j_2,n}}\equiv 1 \ \ {and } \ \  \lim_{n\to \infty} \frac{|x_{j_1,n}- x_{j_2,n}|}{\lambda_{j_1,n}} = \infty.
\end{align}
In this section we always assume that  $u_0 \in H^{d/2-1}$ and ${\rm div}\, u_0=0$.  Let the solution of NS blow up at $T>0$. According to the results in \cite{DoDu09, EsSeSv03, KeKo11, Se12}, we see that $\lim\sup_{t\to T} \|NS(u_0)(t)\|_{\dot H^{d/2-1}} =\infty$. Taking $\tau_n \nearrow T$ so that
\begin{align} \label{sequenceofblowup}
\lim_{n\to \infty} \|NS(u_0)(\tau_n)\|_{\dot H^{d/2-1}} =\infty.
\end{align}
 For convenience, we denote for any $\lambda_{j,n}>0, \ x_{j,n} \in \mathbb{R}^d$,
\begin{align} \label{dilation}
{\Lambda}^{\alpha}_{j,n} f(x):= \frac{1}{\lambda^{\alpha}_{j,n}} f \left( \frac{x- x_{j,n}}{\lambda_{j,n}} \right), \ \ {\Lambda}_{j,n} f ={\Lambda}^{1}_{j,n} f(x).
\end{align}

Let us recall the profile decomposition for a bounded sequence in $\dot H^{d/2-1}$ (see \cite{Ge98}).

\begin{thm} \label{profile}
Let $u_{0,n}$ be a bounded sequence of divergence-free vector fields in $\dot H^{d/2-1}$. Then there exists $\{x_{j,n}\} \subset \mathbb{R}^d$ and $\{\lambda_{j,n}\}\subset (0,\infty)$ which are orthogonal in the sense of \eqref{orthogonal},
and a sequence of divergence-free vector fields   $\{\phi_j\} \subset \dot H^{d/2-1}$ such that
\begin{align} \label{orthogonal decomposition}
u_{0,n} (x) = \sum^J_{j=1} \Lambda_{j,n} \phi_j (x) +  \omega_n^J (x),
\end{align}
where $\omega^J_n $ is a reminder in the sense that
\begin{align} \label{reminder}
\lim_{J\to \infty} \limsup_{n\to \infty}  \|\omega_n^J (x) \|_d =0.
\end{align}
Moreover, we have
\begin{align} \label{orthogonal norm}
\|u_{0,n} (x)\|^2_{\dot H^{d/2-1}}  =  \sum^J_{j=1}  \| \phi_j \|^2_{\dot H^{d/2-1} }  + \|\omega^J_n\|^2_{\dot H^{d/2-1} }   + \epsilon^J_n, \ \ \lim_{n\to\infty}\epsilon^J_n=0, \ \forall \ J.
\end{align}
\end{thm}
Using the conservation law \eqref{conservation} in $L^2$, we can further show that

\begin{lem} \label{weak convergence} {\rm (Weak convergence)}
Let $\{NS(u_0)(\tau_n)\}$ be the solution sequence satisfying \eqref{sequenceofblowup},  $v(\tau_n)= NS(u_0)(\tau_n)/ \|NS(u_0) (\tau_n)\|_{\dot H^{d/2-1}}$. Then $v(\tau_n)$ is weakly convergent to zero in $ H^{d/2-1}$, $L^d$.
\end{lem}

{\bf Proof.}  In view of  \eqref{conservation} and $\lim_{n\to \infty}\|NS(u_0)(\tau_n)\|_{\dot H^{d/2-1}} \to \infty$, we see that $\|v(\tau_n)\|_2 \to 0$.    This implies the weak convergent limit of $\{v(\tau_n)\}$ is zero in $ H^{d/2-1}$, $L^d$. $\hfill\Box$\\

Lemma \ref{weak convergence} indicates that, if we consider the profile decomposition of $v(\tau_n)$, it weakly vanishes in $L^d(\mathbb{R}^d)$ as $\tau_n \nearrow T$.  In view of Theorem \ref{profile}, $v(\tau_n)$ has a profile decomposition in $\dot H^{d/2-1}$:
\begin{align} \label{profile1}
v(\tau_n,x)   = \sum^J_{j=1} \Lambda_{j,n} \phi_j (x) +  \omega_n^J (x),
\end{align}
where $\omega^J_n$ satisfies \eqref{reminder} and \eqref{orthogonal norm}.  Let $j_0 \in \mathbb{N}$. For any $J\geq j_0$, one has \eqref{profile1}. As
$\phi_{j_0} \in \dot H^{d/2 -1}(\mathbb{R}^d) \hookrightarrow L^d(\mathbb{R}^d) $,  H\"older's inequality gives
\begin{align} \label{profileholder}
   \|f \|_{L^p(|\cdot|\leq R_0)} \leq R_0^{d(1/p-1/q)} \|f\|_{L^q (|\cdot|\leq R_0)},\ \ q\geq p,
\end{align}
So we can take a $R_0>0$ such that
\begin{align} \label{profile2}
   \| \phi_{j_0} \|_{L^2(|\cdot|\leq R_0)}   = \frac{1}{2}  \min\{\| \phi_{j_0} \|_{2},  \  1\}.
\end{align}
Since $\lim_{J\to \infty} \lim\sup_{n\to \infty}  \|\omega_n^J (x) \|_d =0,$ we can choose sufficiently large $J_0, N_0 \in \mathbb{N}$ such that

\begin{align} \label{profile3}
 \|\omega_n^J (x) \|_d  < \varepsilon \leq  \frac{ \| \phi_{j_0} \|_{L^2(|\cdot|\leq R_0)}}{100 R_0^{d/2-1}}, \ \ \forall J\geq J_0, \ n\geq N_0.
\end{align}

\begin{lem} \label{lambdavanish} {\rm (Vanishing $\lambda_{j,n}$)}
Let $\{NS(u_0)(\tau_n)\}$ be the solution sequence satisfying \eqref{sequenceofblowup},  $v(\tau_n)= NS(u_0)(\tau_n)/ \|NS(u_0)(\tau_n)\|_{\dot H^{d/2-1}}$.  Let $v(\tau_n)$ has the profile decomposition as in \eqref{profile1} with $J\geq  j_0$. Then we have
\begin{align} \label{profile3a}
 \lambda_{j_0,n}^{d/2-1}   \leq  \frac{4 \|u_0\|_2}{  \min\{\| \phi_{j_0} \|_{2},  \  1\} \|NS(u_0)(\tau_n)\|_{\dot H^{d/2-1}}}, \ \ n\gg 1.
\end{align}
In particular,  we have $\lambda_{j_0,n}\to 0$ as $n\to \infty.$
\end{lem}
{\bf Proof.} We can assume that the reminder term $\omega^J_n$ satisfying \eqref{profile3}. Otherwise we can choose another profile decomposition with $J\geq J_0$.    Let us write $\varphi_{j_0} = \phi_{j_0}  \chi_{\{|\cdot|\leq R_0\}}$.
Observing that
\begin{align} \label{profile sepration}
   \langle \Lambda_{j,n}f, \ \Lambda_{j_0,n} g \rangle = \frac{\lambda^{d-1}_{j_0,n}}{\lambda_{j,n}} \int g(y) f\left( \frac{\lambda_{j_0,n}}{\lambda_{j,n}} y + \frac{x_{j_0,n}- x_{j,n}}{\lambda_{j,n}} \right) dy.
\end{align}
Taking the inner product of \eqref{profile1} and $\varphi_{j_0}$, we have
\begin{align} \label{profile4}
 \langle v(\tau_n), \ \Lambda_{j_0,n} \varphi_{j_0}  \rangle  = \sum^J_{j=1}   \langle \Lambda_{j,n}\phi_j, \ \Lambda_{j_0,n} \varphi_{j_0}  \rangle  +   \langle \omega^J_n, \ \Lambda_{j_0,n} \varphi_{j_0}  \rangle.
\end{align}
In view of \eqref{profile sepration}, we have
\begin{align} \label{profile5}
 \langle  \Lambda_{j_0,n} \phi_{j_0} , \ \Lambda_{j_0,n} \varphi_{j_0}  \rangle  =  \lambda^{d-2}_{j_0,n} \|\varphi_{j_0}\|^2_2 =  \lambda^{d-2}_{j_0,n} \|\phi_{j_0}\|^2_{L^2(|\cdot|\leq R_0)}.
\end{align}
By H\"older's inequality, \eqref{reminder} and \eqref{profileholder},
\begin{align} \label{profile6}
 |\langle  \omega^J_n , \ \Lambda_{j_0,n} \varphi_{j_0}  \rangle |
  & \leq    \|\omega^J_n \|_d   \|\Lambda_{j_0,n} \varphi_{j_0} \|_{d/(d-1)} \nonumber\\
   & \leq     \varepsilon R_0^{d/2-1} \lambda^{d-2}_{j_0,n} \| \phi_{j_0} \|_{L^2(|\cdot|\leq R_0) }.
\end{align}
Since $C^\infty_0 (\mathbb{R}^d)$ is dense in $L^d$, one can choose  $\varphi_j \in C^\infty_0 (\mathbb{R}^d)$ satisfying
\begin{align} \label{profile7}
 \|\phi_j - \varphi_j \|_d   < \varepsilon/J, \ \ j\neq j_0, \ \ j=1,...,J.
\end{align}
We have
\begin{align} \label{profile8}
    |\langle \Lambda_{j,n}\phi_j, \ \Lambda_{j_0,n} \varphi_{j_0}  \rangle|  \leq  & |\langle \Lambda_{j,n}(\phi_j-\varphi_j), \ \Lambda_{j_0,n} \varphi_{j_0}  \rangle| \nonumber\\
  &  +  |\langle \Lambda_{j,n} \varphi_j, \ \Lambda_{j_0,n} \varphi_{j_0}  \rangle|:= I+II.
\end{align}
First, we estimate $I$. By H\"older's inequality, we have
\begin{align} \label{profile9}
 I \leq  & \|\Lambda_{j,n}(\phi_j-\varphi_j)\|_d \|   \Lambda_{j_0,n} \varphi_{j_0} \|_{d/(d-1)} \nonumber\\
  \leq  & \| \phi_j-\varphi_j \|_d  \lambda^{d-2}_{j_0,n} R_0^{d/2-1} \| \varphi_{j_0} \|_{2} \nonumber\\
  \leq  & \varepsilon  J^{-1} \lambda^{d-2}_{j_0,n} R_0^{d/2-1} \| \phi_{j_0} \|_{L^2(|\cdot|\leq R_0)}.
\end{align}
Next, we estimate $II$ and divide the proof into the following three cases.

{\it Case 1.} $\lambda_{j_0,n}= \lambda_{j,n}$ and $\lim_{n\to \infty} |x_{j_0,n}- x_{j,n}|/\lambda_{j,n}=\infty$.  Since $\varphi_j \in C^\infty_0 (\mathbb{R}^d)$, and one can also approximate $\varphi_{j_0}$ by a $C^\infty_0(\mathbb{R}^d)$ function, in view of \eqref{profile sepration}
we see that $II \to 0$ if $n$ is sufficiently large.

{\it Case 2.} $\lim_{n\to \infty} \lambda_{j,n}/ \lambda_{j_0,n}=0$.  Assume that ${\rm supp} \varphi_j \subset B(0, R_1)$. For convenience, we denote $ a_{j,n} =(x_{j_0,n}- x_{j,n})/\lambda_{j_0,n}$.  We have from \eqref{profile sepration}, H\"older's inequality and \eqref{profileholder} that
\begin{align} \label{profile11}
   II & = \frac{\lambda^{d-1}_{j_0,n}}{\lambda_{j,n}} \left| \int \varphi_{j_0}(y)  \varphi_j \left( \frac{\lambda_{j_0,n}}{\lambda_{j,n}} (y + a_{j,n}) \right) dy \right| \nonumber\\
  &  \leq  \frac{\lambda^{d-1}_{j_0,n}}{\lambda_{j,n}}   \|\varphi_{j_0}\|_{L^{d/(d-1)}(|\cdot + a_{j,n}|\leq \lambda_{j,n} R_1/\lambda_{j_0,n}) } \left\| \varphi_j \left( \frac{\lambda_{j_0,n}}{\lambda_{j,n}} (\cdot + a_{j,n})  \right)  \right\|_{d} \nonumber\\
   &  \leq   \lambda^{d-2}_{j_0,n} R^{d/2-1}_{0} \left(\frac{\lambda_{j,n}R_1}{\lambda_{j_0,n}R_0}\right)^{d/2-1} \|\phi_{j_0}\|_{L^2(|\cdot|\leq R_0) } \left\| \varphi_j \right\|_{d}.
   \end{align}
In view of $\lim_{n\to \infty} \lambda_{j,n} / \lambda_{j_0,n}=0$, we see that
\begin{align} \label{profile12}
   \left(\frac{\lambda_{j,n}R_1}{\lambda_{j_0,n}R_0}\right)^{d/2-1}  \left\| \varphi_j \right\|_{d}  \leq \varepsilon /J, \ \ j=1,...,J, \ j\neq j_0
\end{align}
if $n$ is sufficiently large. Hence, it follows from \eqref{profile11} and \eqref{profile12} that
\begin{align} \label{profile13}
   II \leq  \varepsilon J^{-1} \lambda^{d-2}_{j_0,n} R^{d/2-1}_{0}  \|\phi_{j_0}\|_{L^2(|\cdot|\leq R_0) }, \ \ n\gg 1.
\end{align}

{\it Case 3.} $\lim_{n\to \infty} \lambda_{j_0,n} /\lambda_{j,n} =0$.  Still assume that ${\rm supp} \ \varphi_j \subset B(0, R_1)$ and $ b_{j,n} =(x_{j,n}- x_{j_0,n})/\lambda_{j,n}$. Similarly as in Case 2, we have from \eqref{profile sepration}, H\"older's inequality and \eqref{profileholder}  that
\begin{align} \label{profile14}
   II & = \frac{\lambda^{d-1}_{j,n}}{\lambda_{j_0,n}} \left| \int \varphi_{j}(y)  \varphi_{j_0} \left( \frac{\lambda_{j,n}}{\lambda_{j_0,n}} (y + b_{j,n}) \right) dy \right| \nonumber\\
  &  \leq  \frac{\lambda^{d-1}_{j,n}}{\lambda_{j_0,n}}   \|\varphi_{j}\|_{L^{d}(|\cdot + b_{j,n}|\leq \lambda_{j_0,n} R_0/\lambda_{j,n}) } \left\| \varphi_{j_0} \left( \frac{\lambda_{j,n}}{\lambda_{j_0,n}} ( \cdot + b_{j,n})  \right)  \right\|_{L^{d/(d-1)}} \nonumber\\
  &  \leq   \lambda^{d-2}_{j_0,n} R^{d/2-1}_{0}  \|\phi_{j_0}\|_{L^2(|\cdot|\leq R_0) } \left\| \varphi_j \right\|_{L^{d}(|\cdot+ b_{j,n}|\leq \lambda_{j_0,n}R_0/ \lambda_{j,n} ) }.
   \end{align}
Noticing that for any $j=1,...,J$ and $j\neq j_0$, $\varphi_j\in L^d (\mathbb{R}^d)$, in view of the absolute continuity of the integration, $\lim_{n\to \infty} \lambda_{j_0,n}/ \lambda_{j,n}=0$ implies that
\begin{align} \label{profile15}
   \left\| \varphi_j \right\|_{L^{d}(|\cdot+ b_{j,n}|\leq \lambda_{j_0,n}R_0/ \lambda_{j,n} ) } \leq \varepsilon /J
   \end{align}
if $n$ is sufficiently large. Hence, it follows from \eqref{profile14} and \eqref{profile15} that for $1\leq j \leq J, \ j\neq j_0$,
\begin{align} \label{profile16}
   II \leq  \varepsilon J^{-1} \lambda^{d-2}_{j_0,n} R^{d/2-1}_{0}  \|\phi_{j_0}\|_{L^2(|\cdot|\leq R_0) }, \ \ n\gg 1.
\end{align}

Collecting \eqref{profile8}, \eqref{profile9}, \eqref{profile13} and \eqref{profile16}, one has that,
\begin{align} \label{profile17}
 \sum_{j\neq j_0} |\langle \Lambda_{j,n}\phi_j, \ \Lambda_{j_0,n} \varphi_{j_0}  \rangle|     \leq 2 \varepsilon   \lambda^{d-2}_{j_0,n} R^{d/2-1}_{0}  \|\phi_{j_0}\|_{L^2(|\cdot|\leq R_0) }, \ \ n\gg 1.
\end{align}
Now, using \eqref{profile4}--\eqref{profile6} and \eqref{profile17}, we have
\begin{align} \label{profile18}
 |\langle v(\tau_n), \ \Lambda_{j_0,n} \varphi_{j_0}  \rangle|   &  \geq  \lambda^{d-2}_{j_0,n} \left( \|\varphi_{j_0}\|^2_{L^2(|\cdot|\leq R_0)} -  3 \varepsilon   R^{d/2-1}_{0}  \|\phi_{j_0}\|_{L^2(|\cdot|\leq R_0) } \right) \nonumber\\
 &  \geq \frac{1}{2} \lambda^{d-2}_{j_0,n}   \|\phi_{j_0}\|^2_{L^2(|\cdot|\leq R_0)}.
\end{align}
Using the conservation \eqref{conservation}, we have
\begin{align} \label{profile19}
\|u_0\|_2 \| \Lambda_{j_0,n} \varphi_{j_0}  \|_2  & \geq  \|NS(u_0)(\tau_n)\|_2   \| \Lambda_{j_0,n} \varphi_{j_0}  \|_2  \nonumber\\
 & \geq \frac{1}{2} \lambda^{d-2}_{j_0,n} \|NS(u_0) (\tau_n)\|_{\dot H^{d/2-1}} \|\phi_{j_0}\|^2_{L^2(|\cdot|\leq R_0)}.
\end{align}
Hence,
\begin{align} \label{profile20}
\lambda^{d/2-1}_{j_0,n}     \leq  \frac{4 \|u_0\|_2}{\min\{\|\phi_{j_0}\|_{2}, \ 1\} \|NS(u_0) (\tau_n)\|_{\dot H^{d/2-1}}}, \ \ n\gg 1.
 \end{align}
Since $\lim_{n\to \infty} \|u(\tau_n)\|_{\dot H^{d/2-1}}=\infty$, we easily see that $\lim_{n\to \infty} \lambda_{j_0,n}=0$. $\hfill\Box$

\begin{lem}\label{Blowup profile1} {\rm (Blowup profile)} Let $\{NS(u_0) (\tau_n)\}$ be as in Lemma \ref{lambdavanish}.  There exist $\{\alpha_n\}$, $\{\lambda_j(\tau_n)\}$, $\{\phi_j\}$ and $\{x_{j}(\tau_n)\}$  with   $ \alpha_n \to +\infty$, $\lambda_j(\tau_n) \to 0$, $\phi_j \in \dot H^{d/2-1}$ and $x_{j}(\tau_n)\in \mathbb{R}^d$ such that  $u(\tau_n)$ can be decomposed in the following way:
\begin{align} \label{blowupdecompostion}
NS(u_0) (\tau_n)= \alpha_n \left(\sum^J_{j=1} \frac{1}{\lambda_{j}(\tau_n)} \phi_j \left( \frac{x- x_{j}(\tau_n)}{\lambda_{j}(\tau_n)}  \right) + \omega^J_n \right).
\end{align}
where $\omega^J_n$ satisfies $\lim_{J\to \infty}\lim\sup_{n\to\infty}\|\omega^J_n\|_d=0$ and moreover,    $\alpha_n \sim \|NS(u_0)(\tau_n)\|_{\dot H^{d/2-1}}$, $\lambda_j (\tau_n) \lesssim  \|NS(u_0) (\tau_n)\|^{-2/(d-2)}_{\dot H^{d/2-1}} $ and
\begin{align} \label{concentrationpicture}
\sup_{\xi \in \mathbb{R}^d}  \|NS(u_0) (\tau_n)\|_{L^p(|\cdot-\, \xi| \lesssim   \lambda_{j}(\tau_n))}  \gtrsim   \lambda_{j}(\tau_n)^{d/p-1}   \|NS(u_0) (\tau_n)\|_{\dot H^{d/2-1}}.
\end{align}
\end{lem}
{\bf Proof.}
Let us observe that
\begin{align} \label{profile21}
 |\langle v(\tau_n), \ \Lambda_{j_0,n} \varphi_{j_0}  \rangle| \leq \|v(\tau_n)\|_{L^2(|\cdot-x_{j_0,n}| \leq R_0\lambda_{j_0,n})}  \lambda^{d/2-1}_{j_0,n} \| \varphi_{j_0} \|_2,
\end{align}
we have from \eqref{profile18} that
\begin{align} \label{profile22}
\|NS(u_0)(\tau_n)\|_{L^2(|\cdot-x_{j_0,n}| \leq  R_0\lambda_{j_0,n})}
    \geq \frac{1}{4}    \lambda^{d/2-1}_{j_0,n}  \|\varphi_{j_0}\|_{2} \|NS(u_0)(\tau_n)\|_{\dot H^{d/2-1}}.
\end{align}
Again, in view of \eqref{profileholder}, one has that for any $2\leq p \leq \infty$,
\begin{align} \label{profile23}
\|NS(u_0)(\tau_n)\|_{L^2(|\cdot-x_{j_0,n}| \leq  R_0\lambda_{j_0,n})}
    \leq (R_0 \lambda_{j_0,n})^{d/2-d/p} \|NS(u_0)(\tau_n)\|_{L^p(|\cdot-x_{j_0,n}| \leq  R_0\lambda_{j_0,n})}.
\end{align}
Combining \eqref{profile22} and  \eqref{profile23},
\begin{align} \label{profile24}
  \|NS(u_0)(\tau_n)\|_{L^p(|\cdot-x_{j_0,n}| \leq  R_0\lambda_{j_0,n})} \geq \frac{1}{4}  R_0^{d/p-d/2}   \lambda_{j_0,n}^{d/p-1} \|\varphi_{j_0}\|_{2} \|NS(u_0)(\tau_n)\|_{\dot H^{d/2-1}}.
\end{align}
Taking $x_{j,n}= x_j (\tau_n)$, $\lambda_j(\tau_n)= \lambda_{j,n}$, we see that the result follows from the profile decomposition in Theorem \ref{profile}. $\hfill\Box$

\section{Profile decomposition in $L^p$}

 First, let us recall the following theorem concerning the profile decomposition of bounded sequence in $L^p(\R^d)$.
\begin{thm}[\rm \cite{Ko10}] \label{thm:lppd}
 Let $2 \leq p < q,r \leq \infty $, $s_{r,p}:= d(1/r-1/p)$. Let $\{f_n\}^{\infty}_{n=1}$ be a bounded sequence in $L^p(\R^d)$ ,there exists a sequence of profiles $\{\phi_j\}^{\infty}_{j=1} \subset L^p(\R^d) $ and orthogonal sequences of scales and cores $\{(\lambda_{j,n},x_{j,n})\}_{n=1}^{\infty}$ ($j\in \mathbb{N}$), such that up to an extraction in $n$,
\begin{equation}\label{profequ}
f_n =
\sum^J_{j=1} {\Lambda}^{d/p}_{j,n}   \phi_j  + \psi^J_n \,
\end{equation}
and the following properties hold:
\begin{itemize}
\item[\rm (i)] The remainder term $\psi_n^J $ satisfies the smallness condition:
\begin{equation}\label{orthob}
\lim_{J \to \infty}
\big(\limsup_{n \to \infty}\|\psi_n^J\|_{\dot B^{s_{r,p}}_{r,q}}
\big)=0.
\end{equation}
\item[\rm (ii)] The profiles satisfy the following inequality:
\begin{equation}\label{orthoc}
\sum_{j=1}^{\infty}\|\phi_j\|_{p}^p \lesssim \liminf_{n \to \infty}\|f_n\|_{p}^p
\end{equation}
and for each integer $J$,
\begin{equation}\label{orthod}
\|\psi_n^J\|_{p} \lesssim \|f_n\|_{p} + \circ(1) \quad \textrm{as} \quad n\to\infty\,.
\end{equation}
\end{itemize}
\end{thm}
Using the diagonal method, we can assume that, up to an extraction, $\lim_{n\to \infty}\lambda_{j,n}$ exists for all $j\in \mathbb{N}$ (whose limit can be $+\infty$).
We further denote
\begin{align}
& \mathbb{J}_1  =\{j \in \N: \lim_{n\to\infty}\lambda_{j,n} \neq 0, \infty \},  \ \ \mathbb{J}_1 (J)= \mathbb{J}_1 \cap \{j\leq J\} , \nonumber\\
& \mathbb{J}_0 (J) = \{j \in \N: \lim_{n \to \infty}\lambda_{j,n} = 0, \ j\leq J\}, \nonumber\\
& \mathbb{J}_{\infty}(J) = \{j \in \N: \lim_{n \to \infty} \lambda_{j,n}= \infty, \ j\leq J \}, \nonumber\\
 &  \mathbb{J}_1^c (J) = \mathbb{J}_0 (J) \cup \mathbb{J}_{\infty}(J). \nonumber
\end{align}
For $j \in \mathbb{J}_1$, we will simply call $\lambda_{j,n}$ constant scale. Now $f_n$ can be rewritten as
$$
f_n =
\sum_{j \in \mathbb{J}_1(J)} \Lambda^{d/p}_{j,n} \phi_j  + \sum_{j \in \mathbb{J}_1^c (J)}\Lambda^{d/p}_{j,n}\phi_j +  \psi^J_n.
$$
Denote
$$
g_{\eta}(x):= g(x)\chi_{\{\frac{1}{\eta} \leq |g(x)| \leq \eta \} }, \quad \,  g_{\eta^c}(x):= g(x) - g_{\eta}(x).
$$
Next, we show that, if the scales $\lambda_{j,n}$ are very small or very large, then the corresponding profiles  $\Lambda^{d/p}_{j,n} \phi_j$ are very small in appropriate function spaces.
\begin{prop}\label{small}
Let $p,q,r,s_{r,p}$ are as in Theorem~$\ref{thm:lppd}$, then
\begin{itemize}
\item[\rm (1)] $\displaystyle{\lim_{n \to \infty}\left\|\sum_{j\in \mathbb{J}_0 (J)}\lmds\phi_{j,\eta}\right\|_{p_1} = 0}$  \quad
for any \, $ p_1 < p\,\,\textrm{and fixed}\,\,\, J ,\eta \geq 1$\,.

\item[\rm (2)] $\displaystyle{\lim_{n \to \infty}\left\|\sum_{j\in \mathbb{J}_{\infty} (J)}\lmds\phi_{j,\eta}\right\|_{p_2} = 0}$\quad
for any \, $p_2 >p\,\,\textrm{and fixed}\,\,\, J ,\eta \geq 1$\,.

\item[\rm (3)] $\displaystyle{\lim_{J \to \infty}\lim_{\eta \to \infty}\limsup_{n\to \infty}
\left\|\sum_{j \in \mathbb{J}_1^c (J)}\lmds \phi_{j,\eta^c} + \psi_n^J \right\|_{\dot B^{s_{r,p}}_{r,q}} = 0 }$ .
\end{itemize}
\end{prop}

{\bf{Proof.}}    First, we prove (1).  Let $p > p_1$, $J , \eta$ \,fixed, then
$$
\begin{aligned}
\big\|\sjjz{0} \lmds\phi_{j,\eta} \big \|_{p_1}
& \leq \sjjz{0} \lambda_{j,n}^{d(\frac{1}{p_1} - \frac{1}{p})}\|\phi_{j,\eta}\|_{p_1} \\
& \leq \sjjz{0} \lambda_{j,n}^{d(\frac{1}{p_1} - \frac{1}{p})} \eta^{\frac{p}{p_1}-1}\|\phi_{j}\|_{p}^{\frac{p}{p_1}}
\longrightarrow 0 \qquad as \quad n\to \infty.
\end{aligned}
$$
The proof of (2) is similar to (1), and the details of the proof are omitted.

Finaly, we prove (3). By the assumptions on $p,q,r$ and Proposition~\ref{embedbs}, we can easily see $L^p \hookrightarrow \dot B_{r,q}^{s_{r,p}}$, thus
$$
\begin{aligned}
\big\|\sjjzc \lmds \phi_{j,\eta^c} + \psi_n^J \big\|_{\dot B_{r,q}^{s_{r,p}}} & \lesssim \big\|\sjjzc \lmds\phi_{j,\eta^c}
\big\|_{p} + \|\psi_n^J\|_{\dot B_{r,q}^{s_{r,p}}}  \\
&\leq \sjjzc\|\phi_{j,\eta^c}\|_{p} + \|\psi_n^J\|_{\dot B_{r,q}^{s_{r,p}}}
\end{aligned}
$$
By the Lebesgue dominated convergence theorem and the definition of $\phi_{j,\eta^c}$,  we see that
$$
\lim_{\eta \to \infty}\limsup_{n \to \infty}\sjjzc\|\phi_{j,\eta^c}\|_{p}
= 0.
$$
So, we have the result, as desired. $\hfill\Box$

Now, let us denote
\begin{align} \label{j1lambda}
\lim_{n\to \infty} \lambda_{j,n} = \lambda_j, \ \ \Lambda_j \phi_j =  \frac{1}{\lambda^{3/p}_{j}} \phi_j \left( \frac{\cdot}{\lambda_j }\right), \ \  j\in \mathbb{J}_1.
\end{align}
Proposition~\ref{small} indicates that the scales with indices in $\mathbb{J}_0\cup \mathbb{J}_{\infty}$ may not take the major roles to the Navier-Stokes evolution, as they can be considered to be the error term in some function spaces. In fact, we have the following

\begin{thm}[\rm NS evolution of the bounded $\lp$ initial data]\label{structhm}
Fix $3<p<\infty$, let $(u_{0,n})_{n \geq 1}$ be a bounded sequence of divergence-free vector fields in $L^p(\rt)$, whose profile decomposition is given in the
following form as in $Theorem~\ref{thm:lppd}$:
$$
u_{0,n}(x) = \sum_{j=1}^J \Lambda_{j,n}^{3/p}\phi_j(x) + \psi_n^J(x)\,.
$$
We have the following conclusions.
\begin{itemize}
\item[\rm (i)] The existence time of the Navier-Stokes solution associated with initial data $u_{0,n}$ satisfies
\begin{equation}\label{existtime}
\liminf_{n \geq 1}T (u_{0,n}) \geq \widetilde{T} := \inf_{j \in \mathbb{J}_1} T (\Lambda_j \phi_j)
\end{equation}
and $T (u_{0,n}) = \infty$ if $\mathbb{J}_1$ is empty.

\item[\rm (ii)] There exists a $J_0 \in \N$ and $N(J) \in \N$, up to an extraction on $n$, such that
\begin{equation}\label{nsequ}
NS(u_{0,n}) = \sum_{j \in \mathbb{J}_1(J)}NS(\Lambda^{3/p}_{j,n} \phi_j)  + \etl{(\sum_{j \in \mathbb{J}_1^c(J)}\Lambda_{j,n}^{3/p}\phi_j(x) + \psi_n^J(x) )} +R_n^J
\end{equation}
is well-defined for $J>J_0, \, n>N(J) \,\textrm{and} \,\,\, t<\widetilde{T}$,
moreover, for any $T < \widetilde{T}$,
\begin{equation}\label{smrem}
\lim_{J \to \infty}\limsup_{n\to \infty} \|R_n^J\|_{\X} = 0, \ \   \X : = C([0,T],L^p_x) \cap L^{\frac{5p}{3}}((0,T); L^{\frac{5p}{3}}_x).
\end{equation}
%where $X^p_T := C([0,T],L_x^p) \cap L_t^{\frac{5p}{3}}(0,T; L_x^{\frac{5p}{3}})$

\item[\rm (iii)] The solution satisfies the orthogonality in the sense
\begin{equation}\label{soluorth}
\|NS(u_{0,n})(t)\|_{p}^p = \sum_{j \in \mathbb{J}_1(J)} \|NS(\Lambda_{j} \phi_j)(t)\|_{p}^p + \|\lineapp\|_{p}^p + \gamma_n^J(t)
\end{equation}
with $\lim\limits_{J \to \infty}{\limsup\limits_{n \to \infty}|\gamma_n^J(t)|} = 0 \,\,\,  \textrm{for each time} \,\, 0 < t < \widetilde{T}$.
\end{itemize}
\end{thm}

\begin{rem}\label{rk:4}
$\widetilde{T} $ is attainable, provided $\mathbb{J}_1$ is nonempty, i.e. there exists some $j_0 \in \mathbb{J}_1$ such that $\widetilde{T} = T (\Lambda_{j_0}\phi_{j_0})$. In fact, this is a simple consequence of~{\rm{(\ref{lowerbound})}} and~\rm{(\ref{orthoc})}.
\end{rem}

\begin{rem}\label{rk:5}
By~{\rm{(\ref{orthoc})}}, we know $\lim\limits_{j\to \infty} \|\phi_j\|_{p} =0 $, so for any fixed $T < \widetilde{T}$, there exists an index $J_c=J_c(T)$, such that $T \,\|\phi_j\|^{\sigma_p}_{p} < c_0, \, \,j \geq J_c$, where $c_0$ is given in \eqref{lowerbound} of Theorem \ref{th:1} is sufficiently small, then by the fixed point argument in proving Theorem \ref{th:1}, we see that $\|NS( \phi_j)\|_{\mathfrak{X}_T^p} \leq
2C \|\phi_j\|_{p}, \,\, j\geq J_c $. For the details one can refer to  Appendix A.
\end{rem}

The proof of Theorem~\ref{structhm} will be given in Sections~\ref{pfstructhm} and Section~\ref{pfthfour}.

%%%%%%%%%%%%%%%%%%%%%%%%%%%%%%%%%%%%%%%%%%%%%%%%%%%%%%%%%%%%%%%%%%%%%%%%%%%%%%%%%%%%%
%%%%%%%%%%%%%%%%%%%%%%%%%%%%%%%%%%%%%%%%%%%%%%%%%%%%%%%%%%%%%%%%%%%%%%%%%%%%%%%%%%%%%
%%%%%%%%%%%%%%%%%%%%%%%%%%%%%%%%%%%%%%%%%%%%%%%%%%%%%%%%%%%%%%%%%%%%%%%%%%%%%%%%%%%%%

\section{Proof of Theorem ~\ref{mainthm}\label{pfmainthm}}

Denote
\begin{displaymath}
D_c =\{NS(u_0): \,  T (u_0) < \infty,  \quad \limsup_{t \to T (u_0)}(T (u_0)-t)\|NS(u_0)(t)\|^{\sigma_p}_{L^p_x} = M^{\sigma_p}_c \}.
\end{displaymath}
In this Section, we will  prove Theorem~\ref{mainthm} by using Theorem~\ref{structhm}, whose proof are separated into two steps. First, we show that there is a minimal solution $NS(u_0)\in D_c$. Secondly,  using the smoothing effect of the bilinear term of the Navier-Stokes equation, we finally construct a new solution which belongs to both $L^p$ and $\bph$. The proof follows the same ideas as in Poulon~\cite{Po15}.

\subsection{Existence of minimal solution in subcritical $L^p$ \label{lpcrit}}

\begin{prop}\label{idne}
Let $M^{\sigma_p}_c<\infty$.  Then $D_c$ is nonempty.
\end{prop}
{\bf{Proof.}} By the definition of $\CC$, there
exists a sequence of $(u_{0,n})_{n \geq 1}$, such that $T(u_{0,n}) < \infty$ and
\begin{equation}\label{}
\limsup_{t \to T (u_{0,n})}(T (u_{0,n})-t)\|NS(u_{0,n})\|_{p}^{\sigma_p} < \CC + \frac{1}{2n}.
\end{equation}
Further, one can find a time sequence $t_n\, (\mtn - t_n \xrightarrow{n\to \infty} 0)$ verifying
\begin{equation}\label{inequans}
(\mtn -t)\|NS(u_{0,n})\|_{p}^{\sgp} < \CC + 1/n  \quad \, \textrm{for all} \quad t\geq t_n.
\end{equation}
For simplicity, we write $\mu_n = \mtn - t_n$  and define a rescaled initial data $v_{0,n}$
$$
v_{0,n} = \mu_n^{\frac{1}{2}}NS(u_{0,n})(t_n,\mu_n^{\frac{1}{2}}x).
$$
Obviously from ~(\ref{inequans}), $\|v_{0,n}\|_{p}^{\sgp} = \mu_n\|NS(u_{0,n})(t_n)\|_{p}^{\sgp} < \CC + 1/n $,
by uniqueness, the associated solution of NS is
$$
NS(v_{0,n})(\tau,x) = \mu_n^{\frac{1}{2}}NS(u_{0,n})(t_n + \mu_n\tau,\mu_n^{\frac{1}{2}}x)
$$
which are all defined on time interval $[0,1)$. Again, it follows from \eqref{inequans} that
\begin{align}\label{estscans}
(1-\tau)\|NS(v_{0,n})(\tau)\|_{p}^{\sgp}  \leq \CC + 1/n .
\end{align}
Then the desired result is a direct consequence of the following Lemma~\ref{existlp}.

\begin{lem}\label{existlp}
Let $(\epn)_{n\geq 1}$ be a decreasing and vanishing sequence, $(v_{0,n})_{n \geq 1}$ is a bounded sequence in $L^p(\rt)$ and satisfies
\begin{equation}\label{assump}
T (v_{0,n}) =1 , \qquad  \sup_{0< \tau < 1}(1 - \tau)\|NS(v_{0,n})(\tau)\|_{p}^{\sgp} \leq \CC + \epn
\end{equation}
Then the following conclusions hold (passing to a subsequence if necessary)

\begin{itemize}
\item  \textrm{In the profile decomposition of} $v_{0,n}$, there exists a unique scale profile $\Lambda_{k_0}\phi_{k_0}$ such that $NS(\Lambda_{k_0} \phi_{k_0})$ blows up at time $1$.
\item $NS(\Lambda_{k_0} \phi_{k_0})\in D_c$ and
\begin{equation}\label{supcrit}
\sup_{0<\tau<1}(1-\tau)\|NS(\Lambda_{k_0} \phi_{k_0})(\tau)\|_{p}^{\sgp} =
\limsup_{\tau \to 1}(1-\tau)\|NS(\Lambda_{k_0} \phi_{k_0})(\tau)\|_{p}^{\sgp} = \CC .
\end{equation}
\end{itemize}
\end{lem}
{\bf{Proof.}} Assume the profile decomposition of $(v_{0,n})$ has the following form
\begin{equation}\label{}
v_{0,n}(x) = \sjjz{1}  \Lambda^{3/p}_{j,n} \phi_j   + \sjjzc\lmdst\phi_j + \psi_n^J.
\end{equation}
 Recall that $\widetilde T = \inf\limits_{j \in \mathbb{J}_1}T (\Lambda_{j} \phi_j)$, by Remark~\ref{rk:4}, we assume $\widetilde T = T (\Lambda_{k_0} \phi_{k_0}) $ and divide the proof into three steps:

$\rm{(i)}$  We show that $\widetilde T = T (\Lambda_{k_0} \phi_{k_0}) =1 $.
In view of~(\ref{existtime}) in Theorem~\ref{structhm},  one has that $\widetilde T \leq 1$.
Assume $\widetilde T < 1$, for any  $\tau < \widetilde T$,  by (\ref{soluorth}),
$$
\|NS(v_{0,n})(\tau)\|_{p}^p \geq \|NS(\Lambda_{k_0} \phi_{k_0})(\tau)\|_{p}^p - |\gamma_n^J(\tau)|.
$$
Multiplying $(T (\Lambda_{k_0}\phi_{k_0})-\tau)^{\frac{p}{\sgp}}$ on both sides, and using the assumption (\ref{assump}), one has that
\begin{equation}\label{inequ}
\frac{(T (\Lambda_{k_0} \phi_{k_0}) - \tau)^{\frac{p}{\sgp}}}{(1 - \tau)^{\frac{p}{\sgp}}}(\CC + \epn)^{\frac{p}{\sgp}}
\geq \big\{(T (\Lambda_{k_0} \phi_{k_0}) - \tau)^{\frac{1}{\sgp}}\|NS(\Lambda_{k_0} \phi_{k_0})(\tau)\|_{p}\big\}^p  -  |\gamma_n^J(\tau)|.
\end{equation}
Taking limit with regard to $n$ and $J$ on both sides, we have
\begin{equation}\label{inequality}
\frac{(T (\Lambda_{k_0} \phi_{k_0}) - \tau)^{\frac{p}{\sgp}}}{(1 - \tau)^{\frac{p}{\sgp}}}M_c^p
\geq \big\{(T (\Lambda_{k_0} \phi_{k_0}) - \tau)^{\frac{1}{\sgp}}\|NS(\Lambda_{k_0} \phi_{k_0})(\tau)\|_{p}\big\}^p.
\end{equation}
Noting that $\tau < \widetilde T = T (\Lambda_{k_0} \phi_{k_0})$, so we have
\begin{equation*}
 \frac{(\widetilde T - \tau)^{\frac{p}{\sgp}}}{(1 - \widetilde T)^{\frac{p}{\sgp}}} M_c^p \geq \big\{(T (\Lambda_{k_0} \phi_{k_0}) - \tau)^{\frac{1}{\sgp}}\|NS(\Lambda_{k_0} \phi_{k_0})(\tau)\|_{p}\big\}^p.
 \end{equation*}
The left hand side of (\ref{inequality}) can be arbitrarily small if we take $\tau$ sufficently close to $\widetilde T$, however, the right hand side of (\ref{inequality}) has a lower bound as in Theorem~\ref{th:1}, which leads to a contradiction. Hence $\widetilde T = T (\Lambda_{k_0}\phi_{k_0}) =1$.

$\rm{(ii)}$  We prove that $NS(\Lambda_{k_0} \phi_{k_0})\in D_c$.
It follows from the definition of $\CC$ and (\ref{inequality}) that
$$
\limsup_{\tau \to 1}(1-\tau)\|NS(\Lambda_{k_0} \phi_{k_0})(\tau)\|_{p}^{\sgp} =
\sup_{0<\tau<1}(1-\tau)\|NS(\Lambda_{k_0} \phi_{k_0})(\tau)\|_{p}^{\sgp} = \CC.
$$

$\rm{(iii)}$  Profile with constant scale blowing up at time $1$ is unique.  We also prove it by contradiction. Assume there is another constant scale profile $\Lambda_{j_0}\phi_{j_0}$ which blows up at time $1$. Using (\ref{soluorth}) again, we find
$$
\|NS(v_{0,n})(\tau)\|_{p}^p \geq \|NS(\Lambda_{k_0}\phi_{k_0})(\tau)\|_{p}^p + \|NS(\Lambda_{j_0}\phi_{j_0})(\tau)\|_{p}^p - |\gamma_n^J(\tau)|.
$$
Repeating the the argument in $\rm{(i)}$, we can deduce that
$$
M_c^p \geq \big\{(1 - \tau)^{\frac{1}{\sgp}}\|NS(\Lambda_{k_0}\phi_{k_0})(\tau)\|_{p}\big\}^p + \big\{(1 - \tau)^{\frac{1}{\sgp}}\|NS(\Lambda_{j_0}\phi_{j_0})(\tau)\|_{p}\big\}^p.
$$
Due to $NS(\Lambda_{k_0}\phi_{k_0}) \in D_c$, we must have $NS(\Lambda_{j_0}\phi_{j_0}) = 0$, so $\phi_{j_0}= 0$. The proof of Lemma~\ref{existlp} is finished. $\hfill\Box$

\subsection{Existence of solutions in $D_c\cap \bph$\label{besovcrit}}

In general, we could not expect that the linear term $ \etl u_0$ belongs to $\bph$, provided we only assume the initial data $u_0 \in L^p$. To get the desired result, we need to use the regularization effect of the biliner term $\mathscr{B}(u,u)$. In fact, based on our assumption, the following proposition holds.

\begin{prop}\label{biregu}
Let $u_0$ be in $L^p$, where $3 <p< \infty$, the associated Navier-Stokes solution $u(t):=NS(u_0)(t)$ satisfies a minimal blow up condition, i.e.
$$
\etu < \infty, \quad M^{\sgp} := \sup_{0<t <\etu} (\etu - t)\|NS(u_0)(t)\|_{p}^{\sgp} < \infty , \quad  \sgp = \frac{2}{1-3/p}
$$
then we have
\begin{equation}\label{inequabs}
 \sup_{0<t < \etu}\|\mathscr{B}(u,u)(t)\|_{\bph} < \infty
\end{equation}
Furthermore, if \, $0 < s < 1- 3/p$, then
\begin{equation}\label{inequanbs}
\sup_{0<t < \etu}(\etu - t)^{\frac{1}{2}(s +1 - 3/p)} \|\mathscr{B}(u,u)(t)\|_{\dot B^s_{p,\infty}}< \infty
\end{equation}
\end{prop}
{\bf{Proof.}} First, we prove (\ref{inequabs}). It follows from \eqref{exp-decay 0} that
\begin{equation}\label{estimone}
\begin{aligned}
\|\triangle_j \mathscr{B}(u,u)\|_{L^{p/2}} & \lesssim \int_0^t 2^j e^{-c(t-s)2^{2j}} \|\triangle_j(u \otimes u)\|_{L^{p/2}}\,ds \\
& \lesssim \int_0^t  2^j e^{-c(t-s)2^{2j}} \|u(s) \|^2_{p}\,ds
\end{aligned}
\end{equation}
Since $u(t)$ is a blowup solution of type-I,
$$
\|u(t)\|_{p}^2  \leq \frac{M^2}{(\etu - t)^{1-3/p}}
$$
Hence
\begin{equation}\label{estimtwo}
\begin{aligned}
2^{j(-1 + \frac{6}{p})}\|\triangle_j \mathscr{B}(u,u)\|_{L^{p/2}} &\lesssim
\int_0^t 2^{j\frac{6}{p}} e^{-c(t-s)2^{2j}} \frac{M^2}{(\etu - s)^{1-3/p}} \,ds.
\end{aligned}
\end{equation}
Let us observe that
$$
\begin{aligned}
 2^{j\frac{6}{p}} \int_0^t \chi_{\{\etu -s \leq 2^{-2j}\}}e^{-c(t-s)2^{2j}} \frac{M^2}{(\etu - s)^{1-3/p}} \,ds \lesssim M^2
\end{aligned}
$$
and
$$
\begin{aligned}
\int_0^t \chi_{\{\etu -s \geq 2^{-2j}\}}e^{-c(t-s)2^{2j}} \frac{M^2}{(\etu - s)^{1-3/p}} \,ds \\
 \lesssim 2^{2j} M^2\int_0^t e^{-c(t-s)2^{2j}} \,ds
  \lesssim M^2.
\end{aligned}
$$
(\ref{inequabs}) is obtained.

Next, we prove (\ref{inequanbs}).  We divide the proof into two cases.

Case 1.   $\etu - t < 2^{-2j}$.  By hypothesis, we see that
$$
\begin{aligned}
(\etu - t)^{\frac{1}{\sgp}}\|\mathscr{B}(u,u)(t)\|_{p} &\lesssim (\etu - t)^{\frac{1}{\sgp}}\|u(t)\|_{p} + (\etu - t)^{\frac{1}{\sgp}}\|\etl u_0\|_{p} \\
& \le M + \etu^{\frac{1}{\sgp}} \|u_0\|_{p} < \infty.
\end{aligned}
$$
For $s > 0$, note $\sgp = \frac{2}{1-3/p}$  and $2^j < (\etu - t)^{-\frac{1}{2}}$
$$
\begin{aligned}
2^{js}\|\triangle_j \mathscr{B}(u,u)(t)\|_{p}
 \lesssim
 \frac{1}{(\etu - t)^{\frac{1}{2}(s+ 1-3/p)}}
\end{aligned}
$$
Case 2. $\etu - t \geq 2^{-2j}$.   By Bernstein's inequality and the assumption,
$$
\|\triangle_j \mathscr{B}(u,u)(t)\|_{p} \lesssim \int_0^t 2^{j(1+3/p)} e^{-c(t-\tau)2^{2j}} \frac{M^2}{(\etu - \tau)^{1-3/p }} \,\mathrm{d\tau}.
$$
Multiplying it by $2^{js}$ and using the fact that $\tau < t$ and $s -1 + 3/p < 0$,
$$
\begin{aligned}
2^{js}\|\triangle_j \mathscr{B}(u,u)(t)\|_{p} &\lesssim \int_0^t 2^{j(1+s+3/p)} e^{-c(t-\tau)2^{2j}} \frac{M^2}{(\etu - t)^{1-3/p }} \,\mathrm{d\tau}\\
& \lesssim \frac{M^2}{(\etu - t)^{1-3/p }} 2^{j(s-1 + 3/p)} \\
& \leq \frac{1}{(\etu -t)^{\frac{1}{2}(s + 1 - 3/p)}}
\end{aligned}
$$
Combining $(1) \,\textrm{and}\, (2)$, we finally obtain (\ref{inequanbs}). $\hfill\Box$\\

Next we prove a lemma which  allows us to exclude profiles without constant scales in the profile decomposition of a bounded sequence in $\lp(\rt)$.%carry out the same procedure as E.Poulon~\cite{poulon}.
\begin{lem}\label{Lem:prd} Let $3< p < \infty, 0<s < 1-3/p$, then the following statements hold:
\begin{itemize}
\item [\rm (1)]
 Assume $(f_n)_{n \geq 1}$ is a bounded sequence in $\bph \cap\, \lp$, s.t. $ {\limsup\limits_{n\to \infty} \|f_n\|_{\dot B_{p,\infty}^0}> 0 }$, then  there are no scales which tend to infinity in the $\lp$ profile decomposition of $f_n$.
\item [\rm (2)]
 Assume $(f_n)_{n \geq 1}$ is a bounded sequence in $\dot B^s_{p,\infty} \cap\, \lp$, s.t. $ \limsup\limits_{n\to \infty} \|f_n\|_{\dot B_{p,\infty}^0} > 0 $, then  there are no scales which tend to zero in the $\lp$ profile decomposition of $f_n$.
\end{itemize}
\end{lem}
{\bf{Proof :}} First, we argue by contradiction to prove (1). Assume for a contrary that there exists a scale $\lambda_{k,n}\xrightarrow{n\to\infty} \infty$.
By Theorem~\ref{thm:lppd}, let the profile decomposition of $f_n$ has the following form
$$
f_n = \sum_{j=1}^J \frac{1}{\lambda_{j,n}^{3/p}}\phi_j \big (\frac{x-x_{j,n}}{\lambda_{j,n}}\big ) + \psi_n^J.
$$
It follows from $\displaystyle \limsup_{n\to \infty} \|f_n\|_{\dot B_{p,\infty}^0} > 0 $  that $\displaystyle f_n \not\equiv 0$,
 As we know, each $\phi_j$ is a weak limit of some rigid
tranformation of $f_n$, more precisely, we have
$$
\lambda_{j,n}^{3/p}f_n(\lambda_{j,n}x + x_{j,n}) \rightharpoonup \phi_j
$$
in the sense of tempered distribution. Without loss of generality, we assume $\phi_j \neq 0 \,\, \textrm{for each}\,\, j$. Observing
$$
\lambda_{k,n}^{3/p}\|f_n(\lambda_{k,n}x + x_{k,n})\|_{\bph} \equiv \lambda_{k,n}^{3/p-1}\|f_n\|_{\bph} \xrightarrow{n\to \infty} 0
$$
due to the hypothesis on $f_n$ and $\lambda_{k,n}\xrightarrow{n\to\infty} \infty$. By lower semi-continuity, we get
$$
\|\phi_k\|_{\bph} \leq \liminf_{n\to \infty} \lambda_{k,n}^{3/p}\|f_n(\lambda_{k,n}x + x_{k,n})\|_{\bph} = 0
$$
so $\phi_k = 0 $, this contradicts our assumption.  The proof of (2) proceeds in a similar way as that of (1) and the details are omitted. $\hfill\Box$\\

\noindent {\bf Proof of Theorem \ref{mainthm}.} We only outline the main procedures in the proof of existence of a critical element in $\bph$, as the argument is quite similar to  ~\cite{Po15}. Indeed, we already have a minimal solution  $\Psi = NS(\Psi_0)\in D_c$  from Proposition~\ref{idne},  so there exists a time sequence
$t_n \to T (\Psi_0)$ and $\epn \to 0$ such that
$$
\lim_{n \to \infty}(T (\Psi_0) - t_n)\|NS(\Psi_0)(t_n)\|_{p}^{\sgp} = \CC
$$
and
\begin{equation}\label{crittime}
(T (\Psi_0) - t)\|NS(\Psi_0)(t)\|_{p}^{\sgp} \leq \CC + \epn   \quad \textrm{for any}\quad t \geq t_n
\end{equation}
as before,we define a rescaled initial data $v_{0,n}$
$$
v_{0,n}(x) = (T (\Psi_0) - t_n)^{\frac{1}{2}}NS(\Psi_0)(t_n,(T (\Psi_0) - t_n)^{\frac{1}{2}}x)
$$
denote $\tau_n = T (\Psi_0) - t_n$ for short,   we get
$$
v_{0,n}(x) = \tau_n^{\frac{1}{2}} e^{t_n \Delta}\Psi_0(\tau_n^{\frac{1}{2}}x) + \tau_n^{\frac{1}{2}}\mathscr{B}(\Psi_0,\Psi_0)(t_n,\tau_n^{\frac{1}{2}}x)
$$
Due to $\tau_n \to 0$,~(\ref{crittime})~and Proposition~\ref{biregu}, we can easily get
$$
\begin{aligned}
\tau_n^{\frac{1}{2}}\| e^{t_n \Delta}\Psi_0(\tau_n^{\frac{1}{2}}x)\|_{p} \to 0,  \qquad n\to \infty  \\
\rho_n := \tau_n^{\frac{1}{2}}\mathscr{B}(\Psi_0,\Psi_0)(t_n,\tau_n^{\frac{1}{2}}x) \in \lp \cap \bph \cap \dot B^s_{p,\infty}
\end{aligned}
$$
 Next we take $\lp$ profile decomposition of $\rho_n$,
$$
\rho_n = \sum_{j=1}^J \lmdst\phi_j(x) + \psi_n^J
$$
then come back to $v_{0,n}$
$$
v_{0,n}(x) = \sum_{j=1}^J \lmdst\phi_j(x) + \psi_n^J +
 \tau_n^{\frac{1}{2}} e^{t_n \Delta}\Psi_0(\tau_n^{\frac{1}{2}}x)
$$
In this circumstance,
$\displaystyle \psi_n^J +
 \tau_n^{\frac{1}{2}} e^{t_n \Delta}\Psi_0(\tau_n^{\frac{1}{2}}x)$  can be considered as the remainders, and according to Lemma~\ref{Lem:prd}, in the profile decomposition of $\rho_n$, only profiles with constant scales are left, and each $\phi_j(x)$ is the weak limit of $\lambda^{3/p}_{j,n}\rho_n(\lambda_{j,n}x +x_{j,n})$, so all these profiles lies in $ \lp \cap \bph \cap \dot B^s_{p,\infty}$. Now it can be easily checked that $(v_{0,n})_{n \geq 1}$ satisfy the assumptions of  Lemma~\ref{existlp}, hence there exists a unique scale profile $\Lambda_{l_0}\phi_{l_0} \in \lp \cap \bph$ such that $NS(\Lambda_{l_0}\phi_{l_0}) \in D_c$  and blows up at time $1$. Moreover for any $0<\tau < 1$, we have
 $$
 \begin{aligned}
 \sup_{0<\tau <1} \|NS(\Lambda_{l_0}\phi_{l_0})(\tau)\|_{\bph} &= \sup_{0<\tau <1} \|e^{\tau\Delta} \Lambda_{l_0}\phi_{l_0} + \mathscr{B}(\Lambda_{l_0}\phi_{l_0},\Lambda_{l_0}\phi_{l_0})(\tau)\|_{\bph} \\
 & \lesssim \|\phi_{l_0}\|_{\bph} + \sup_{0<\tau <1} \| \mathscr{B}(\Lambda_{l_0}\phi_{l_0},\Lambda_{l_0}\phi_{l_0})(\tau)\|_{\bph} < \infty
 \end{aligned}
 $$
 where we have used Proposition~\ref{biregu}. Therefore we finish the proof of Theorem~\ref{mainthm}. $\hfill\Box$

%%%%%%%%%%%%%%%%%%%%%%%%%%%%%%%%%%%%%%%%%%%%%%%%%%%%%%%%%%%%%%%%%%%%%%%%%%%%%%%%%%%%%
%%%%%%%%%%%%%%%%%%%%%%%%%%%%%%%%%%%%%%%%%%%%%%%%%%%%%%%%%%%%%%%%%%%%%%%%%%%%%%%%%%%%%
%%%%%%%%%%%%%%%%%%%%%%%%%%%%%%%%%%%%%%%%%%%%%%%%%%%%%%%%%%%%%%%%%%%%%%%%%%%%%%%%%%%%%
\section{Proof of Theorem~\ref{structhm}\label{pfstructhm}}
This section is devoted to the proof of $\mathrm{(i)}\, \textrm{and} \,\mathrm{(ii)} $  in Theorem~\ref{structhm}. Let $(u_{0,n})_{n \geq 1}$ be a bounded sequence in $\lp$ and the corresponding profile decomposition is given below
\begin{equation}\label{sfinitial}
u_{0,n} = \sum_{j \in \mathbb{J}_1 (J)}\Lambda^{3/p}_{j,n}  \phi_j  + \sum_{j \in \mathbb{J}_1^c(J)}\Lambda^{3/p}_{j,n}\phi_j +  \psi^J_n \,,
\end{equation}
we look for an approximative solution to the genuine solution $NS(u_{0,n})$, let
$$
NS(u_{0,n}) = \U + R_n^J
$$
where we denote
$$
\U = \sum_{j \in \mathbb{J}_1 (J)}NS(\Lambda^{3/p}_{j,n}\phi_j) + \lineapp
$$
obviously, $\U$ is defined on $\displaystyle [0,\widetilde T_1),\,\widetilde T_1:= \inf_{j\in \mathbb{J}_1} T (\Lambda^{3/p}_{j,n}\phi_j)$.
$R_n^J$ is a divergence free vector field, which is defined on $[0, T_n),\,T_n := \min\{\widetilde T_1, \mtn\}$ and satisfies the perturbation equation below
\begin{equation}\label{pertureq}\displaystyle
\left\{\begin{array}{ll}  \displaystyle
&\partial_t \re
+ \re \cdot \nabla \re  - \Delta \re + Q(\re,\U) =   -\nabla \cdot G_n^J - \nabla p_n^J \\
 &\mbox{div}\, \re = 0 \\
& R_{n |t=0}^J   = 0
\end{array}\right.
\end{equation}
where $Q(a,b):= a\cdot \nabla b + b \cdot \nabla a =2 \nabla \cdot(a\otimes b)$ for any divergence free vector fields $a,\,b$  and
\begin{equation}\label{forterm}
\begin{aligned}
G_n^J & = \sum_{j,k \in \mathbb{J}_1(J):j \neq k} \sol{j}\otimes \sol{k}  \\
& +\quad  2\sjjz{1} \sol{j}\otimes \lineapp \\
& +\quad \lineapp \otimes \lineapp
\end{aligned}
\end{equation}

Now $\mathrm{(i)}\, \textrm{and} \,\mathrm{(ii)} $ of Theorem~\ref{structhm} will be consequences of the following two lemmas
and Proposition~\ref{perturprop} in the appendix. We remark that by the lower semi-continuity in Corollary~\ref{lscontinuous}, for arbitrary
$\varepsilon >0$, there exists a $N:=N(\varepsilon)$, such that $n > N$, $\widetilde T_1:= \inf_{j\in \mathbb{J}_1} T(\Lambda^{3/p}_{j,n} \phi_j)
\geq \wt:= \inf_{j \in \mathbb{J}_1} T(\Lambda_j \phi_j) - \varepsilon $.

\begin{lem}\label{bddri}
For any $\displaystyle 0< T < \wt := \inf_{j\in \mathbb{J}_1 } T (\Lambda_j\phi_j) $, there exists constant $C>0$ such that, up to an extraction to $n\in \mathbb{N}$,
$$
\lim_{J\to \infty}\limsup_{n\to \infty}\|\U\|_{\X} \leq C.
$$
\end{lem}

\begin{lem}\label{smfor}
Let $3 < p< \infty$, $ 1/\beta(p)  =  1/2  -  3/10p $, $T$ be as in Lemma \ref{bddri}. Then we have
$$
\lim_{J\to \infty}\limsup_{n\to \infty}\|G_n^J\|_{L_T^{\be}L_x^{\al}}  = 0.
$$
\end{lem}

Let us turn to the proof of Lemma~\ref{bddri} and Lemma~\ref{smfor}. First, we give a lemma, which will be used in the sequel.
\begin{lem}\label{Lem:orth}
Let $J \in \N, \,3<p<\infty$,
$r=5p/3$. Let $\phi_j  \in L^p$ and $(x_{j,n},\lambda_{j,n})_{n\geq 1}$ be orthogonal scales and cores in the sense of ~$\mathrm{(\ref{orthoa})}$, then for any $
T \leq \inf_{j\in \mathbb{J}_1(J)}  T_j$, $T_j < T(\Lambda_j \phi_j)$ is arbitrary. Then there exists $N(T, J)$ such that for any $n> N(T, J)$,
\begin{equation}\label{inequorth}
\left\| \sum_{j\in \mathbb{J}_1(J) } NS (\lmdst \phi_j) \right\|^r_{L_t^r ((0,T),L_x^r(\rt))} \leq \sum_{j\in \mathbb{J}_1(J) }    \big\| NS (\Lambda_j \phi_j) \big\|^r_{L^r_t((0,T_j),L_x^r(\rt))} + \varepsilon(J,n),
\end{equation}
where $\displaystyle \lim_{n\to \infty}\varepsilon(J,n) = 0$.
Furthermore, for any $0\leq t \leq T$, up to an extraction to $n\in \mathbb{N}$,
\begin{equation}\label{infinequorth}
\left\| \sum_{j\in \mathbb{J}_1(J) } NS (\lmdst \phi_j)(t)  \right\|^p_{p} = \sum_{j\in \mathbb{J}_1(J) } \| NS (\Lambda_j \phi_j) (t,\cdot)\|^p_{p} + \sigma_{J,n}(t),
\end{equation}
where $\displaystyle \lim_{n\to \infty}\sup_{0<t < T}|\sigma_{J,n}(t)| = 0$.
\end{lem}

Using Lemma \ref{Lem:orth}, we can  prove Lemma~\ref{bddri}. \\

\noindent {\bf{Proof of Lemma~\ref{bddri}.}} By the definition of $\U$,
$$
\|\U\|_{\X} \leq \| \sjjz{1} \sol{j}\|_{\X} + \|\lineapp\|_{\X} :=
\uppercase\expandafter{\romannumeral1} + \uppercase\expandafter{\romannumeral2}
$$
where
$$
\begin{aligned}
  \uppercase\expandafter{\romannumeral1}
&= \| \sjjz{1} \sol{j}\|_{L_T^{\infty}L_x^p} +
\| \sjjz{1} \sol{j}\|_{L_T^{\frac{5p}{3}}L_x^{\frac{5p}{3}}}.
\end{aligned}
$$
From~(\ref{infinequorth}) it follows that
$$
\| \sjjz{1} \sol{j}\|^p_{L_T^{\infty}L_x^p} \leq \sjjz{1}\|NS(\Lambda_j \phi_j)\|^p_{L_T^{\infty}L_x^p} + \|\sigma_{J,n}(t)\|_{L_T^{\infty}}.
$$
Taking $n \to \infty$,  up to an extraction to $n$, one has that
\begin{equation}\label{inequasum}
\limsup_{n \to \infty}\| \sjjz{1} \sol{j}\|^p_{L_T^{\infty}L_x^p} \leq \sjjz{1}\| NS(\Lambda_j \phi_j) \|^p_{L_T^{\infty}L_x^p}.
\end{equation}
In view of Remark~\ref{rk:5}, we see that there exists a $J_c\in \mathbb{N}$  such that for $j \geq J_c$,
$$
\|NS(\Lambda_j \phi_j)\|_{{L_T^{\infty}L_x^p}}   \leq 2 C\|\phi_j\|_{p}.
$$
Therefore, the left hand side of~(\ref{inequasum}) can be bounded by
\begin{equation}\label{conscabd}
\sum_{j \in \mathbb{J}_1 ( J_c) } \|NS(\phi_j)\|^p_{{L_T^{\infty}L_x^p}} + C \sum_{j=1}^{\infty} \|\phi_j\|^p_{p}
\end{equation}
as indicated by (\ref{orthoc}),\,$\displaystyle \sum_{j=1}^{\infty} \|\phi_j\|^p_{p} \lesssim \liminf_{n\to \infty}\|u_{0,n}\|^p_{p}$ and $T < \widetilde T$,
thus
\begin{equation}\label{constbd}
\lim_{J\to\infty}\limsup_{n \to \infty}\| \sjjz{1} \sol{j}\|^p_{L_T^{\infty}L_x^p} < \infty.
\end{equation}
Using the same way as above,  in view of~(\ref{inequorth})and $l^p\hookrightarrow l^{\frac{5}{3}p} $ we have
\begin{equation}\label{constbda}
\begin{aligned}
\limsup_{n\to \infty}  \| \sjjz{1} \sol{j}\|^{\frac{5p}{3}}_{L_T^{\frac{5p}{3}}L_x^{\frac{5p}{3}}}
 & \lesssim \sum_{j\in \mathbb{J}_1(J_c)} \|NS(\Lambda_j \phi_j)\|^{\frac{5}{3}p}_{L_T^{\frac{5p}{3}}L_x^{\frac{5p}{3}}} + \sum_{j \geq 1}\|\phi_j\|^{\frac{5p}{3}}_{p} < \infty.
\end{aligned}
\end{equation}
Therefore, we get the estimate of $\uppercase\expandafter{\romannumeral1}$. Now we estimate $\uppercase\expandafter{\romannumeral2}$. By the estimate of $e^{t\Delta}$ and Theorem \ref{thm:lppd} we have
$$
\begin{aligned}
\uppercase\expandafter{\romannumeral2}
& \lesssim \|u_{0,n}\|_{p} + \|\sjjz{1}\Lambda^{3/p}_{j,n}\phi_j\|_{p}.
 \end{aligned}
$$
It follows from~(\ref{infinequorth}) that
\begin{equation}\label{linetmbd}
\lim_{J\to\infty}\limsup_{n\to \infty}\|\lineapp\|_{\X} \lesssim \sup_{n \geq 1}\|u_{0,n}\|_{p} < \infty.
\end{equation}
Combining the estimate on $\uppercase\expandafter{\romannumeral1}\, \textrm{and} \,\uppercase\expandafter{\romannumeral2}$,
we finally complete the proof of Lemma~\ref{bddri}. $\hfill\Box$\\

Next, we prove Lemma~\ref{Lem:orth} by following the ideas in Lemma 2.6 of \cite{GaKoPl13}. However, in the subcritical cases $\Lambda^{3/p}_{j,n}$ has no scaling invariance for the solutions of NS, i.e., $NS(\Lambda^{3/p}_{j,n}u_0) \neq \Lambda^{3/p}_{j,n} NS(u_0)$.\\

\noindent {\bf{Proof of Lemma~\ref{Lem:orth}.}} For $r = 5p/3 $,   we have
\begin{align}
\bigg\|\sum_{j\in \mathbb{J}_1(J)} NS(\lmdst \phi_j) \bigg\|^r_{L^r_T L_x^r} &= \bigg\| \big|\sum_{j\in \mathbb{J}_1(J)} NS(\lmdst \phi_j)\big|^r \bigg\|_{L_T^{1}L_x^1}  \nonumber\\
& \leq \sum_{j\in \mathbb{J}_1(J)} \big\| NS(\lmdst \phi_j) \big\|^r_{L^r_T L_x^r} + \varepsilon(J,n), \label{conver0}
\end{align}
where
$$
\begin{aligned}
\varepsilon(J,n): &= \bigg\|\big|\sum_{j\in \mathbb{J}_1(J)} NS(\lmdst \phi_j)\big|^r - \sum_{j\in \mathbb{J}_1(J)} \big| NS(\lmdst \phi_j) \big|^r \bigg\|_{L_T^{1}L_x^1} \\
& \leq C_{J,p} \sum_{  j,k\in \mathbb{J}_1(J), \, j \neq k} \big\| |NS(\lmdst \phi_j)| |NS(\lmdp{k} \phi_k)|^{r-1}\big\|_{L_T^{1}L_x^1}.
\end{aligned}
$$
Here we have used an elementary inequality (cf.,e.g.,\cite{Ge98})
\begin{equation}\label{elemineq}
\left ||\sum_{j=1}^L a_j |^m - \sum_{j=1}^L |a_j|^m \right| \leq C_{L,m}\!\!\sum_{1\leq j,k \leq L:j \neq k}|a_j||a_k|^{m-1},\,\, 1< m < \infty.
 \end{equation}
In order to show the result, we divide the proof into the following three steps.

{\it Step 1.} For any $j\in \mathbb{J}_1$ and $\phi_j \in L^p$, we show that
\begin{align}
\lim_{n\to \infty} \|\lmdst \phi_j -(\Lambda_j \phi_j) (\cdot-x_{j,n})\|_p =0. \label{conver}
\end{align}
For any $\varepsilon>0$, one can choose $\varphi_j \in C^\infty_0 (\mathbb{R}^d)$ satisfies
$$
\|\varphi_j- \phi_j\|_p < \varepsilon.
$$
We have
\begin{align}
 \|\lmdst \phi_j & -  (\Lambda_j \phi_j) (\cdot-x_{j,n})\|_p  \nonumber\\
 & \leq \|\lmdst \varphi_j -  (\Lambda_j \varphi_j) (\cdot-x_{j,n})\|_p + \|\lmdst (\varphi_j -  \phi_j)  \|_p + \|\Lambda_j (\varphi_j -  \phi_j)  \|_p \nonumber\\
  & \leq \| (\lambda_j/\lambda_{j,n})^{3/p} \varphi_j(\lambda_j \cdot /\lambda_{j,n} )  -   \varphi_j \|_p + 2\varepsilon.  \label{conver1}
\end{align}
Noticing that $\lim_{n\to \infty} \lambda_j/\lambda_{j,n} =1 $, from the uniform continuity of $\varphi_j$ we can obtain that
\begin{align}
 \lim_{n\to \infty} \| (\lambda_j/\lambda_{j,n})^{3/p} \varphi_j(\lambda_j \cdot /\lambda_{j,n} )  -   \varphi_j \|_p  =0. \label{conver2}
\end{align}
It follows from \eqref{conver1} and \eqref{conver2} that \eqref{conver} holds true.

{\it Step 2.} For any $j\in \mathbb{J}_1$ and $\phi_j \in L^p$, we show that
\begin{align}
\lim_{n\to \infty} \|NS(\lmdst \phi_j) - NS(\Lambda_j \phi_j) (\cdot-x_{j,n})\|_{\X} =0. \label{1conver}
\end{align}
In fact,  in view of the translation invariance of NS, one can regard $x_{j,n}=0$ in \eqref{1conver}. Now, let $v^{n}:= NS(\lmdst \phi_j) - NS(\Lambda_j \phi_j) (x) $, obviously it satisfies the following perturbation equation
$$
 \quad \left\{
\begin{array}{ll}
\partial_t v^{n} & + (v ^{n}\cdot \nabla) v^{n} -\Delta v^{n} +(v^{n} \cdot \nabla) NS(\Lambda_j\phi_j)+ (NS(\Lambda_j\phi_j) \cdot \nabla) v^{n} =  -\nabla p, \   \mathrm{div} \, v^{n} =  0,
\\
v^{n}_{|t=0} &= \Lambda^{3/p}_{j,n} \phi_j -\Lambda_j \phi_j.
\end{array}
\right.
$$
Applying Proposition \ref{perturprop} and \eqref{conver},  we see that for sufficiently large $n> N(T)$,
$$
\|v^{n}\|_{\X} \lesssim  \|\Lambda^{3/p}_{j,n} \phi_j -\Lambda_j \phi_j\|_p   \exp(C\,T^{ 5p(1-\fp)/6}\|  NS(\Lambda_j\phi_j) \|^{5p/3}_{L_T^{5p/3}L_x^{5p/3}}):= \delta (n, T) \to 0,
$$
Hence, we have \eqref{1conver} and
\begin{align} \label{2conver}
\| NS(\lmdst \phi_j)\|_{\X}  \leq  \| NS(\Lambda_j \phi_j)\|_{\X}  + \delta (n, T).
\end{align}
So, by \eqref{2conver} and~(\ref{constbda}),
\begin{align} \label{3conver}
\sum_{j\in \mathbb{J}_1(J)} \big\| NS(\lmdst \phi_j) \big\|^r_{L^r_T L_x^r}  \leq  \sum_{j\in \mathbb{J}_1(J)} \big\| NS(\Lambda_j \phi_j) \big\|^r_{L^r_T L_x^r} + \rho (n, J, T), \ \ \rho(n, J, T)\to 0.
\end{align}
{\it Step 3.} We show that
\begin{align} \label{0conver}
\lim_{n\to\infty} \big\| |NS(\lmdst \phi_j)| |NS(\lmdp{k} \phi_k)|^{r-1}\big\|_{L_T^{1}L_x^1} =0,\quad  \textrm{for each} \quad j \neq k, \ j,k\in \mathbb{J}_1(J).
\end{align}
By \eqref{1conver}, it suffices to show that
\begin{align} \label{4conver}
\lim_{n\to\infty} \big\| |NS(\Lambda_j \phi_j)(\cdot- x_{j,n})| |NS(\Lambda_k \phi_k) (\cdot- x_{k,n})|^{r-1}\big\|_{L_T^{1}L_x^1} =0,\quad  \textrm{for each} \quad j \neq k.
\end{align}
Since $C^\infty_0((0,T)\times \mathbb{R}^d)$ is dense in $L^r_T L^r_x$, it suffices to show \eqref{4conver} by assuming that $NS(\Lambda_j \phi_j) \in  C^\infty_0((0,T)\times \mathbb{R}^d)$ for all $j\in \mathbb{J}_1(J)$. Noticing that $|x_{j,n}-x_{k,n}| \to \infty$, we see that the left hand side of \eqref{4conver} is zero when $n$ is sufficiently large. Hence, we have \eqref{0conver}.

Now \eqref{inequorth} follows from the estimates of Steps 1--3. Next, we turn to prove~(\ref{infinequorth}).
\begin{align}
\bigg\|\sum_{j\in \mathbb{J}_1(J)} NS(\lmdst \phi_j) \bigg\|^p_{p} := \bigg\|\sum_{j\in \mathbb{J}_1(J)} NS(\Lambda_j \phi_j)(\cdot- x_{j,n}) \bigg\|^p_{p}   + \alpha_1 (t, J,n). \label{conver00}
\end{align}
In view of \eqref{1conver},~(\ref{inequasum})~and~(\ref{conscabd}), we see that
\begin{align}
\lim_{n\to \infty} \sup_{0\le t\le T} |\alpha_1 (t, J,n)| =0.  \label{conver01}
\end{align}
Let us write
\begin{align}
\bigg\|\sum_{j\in \mathbb{J}_1(J)} NS(\Lambda_j \phi_j)(\cdot- x_{j,n}) \bigg\|^p_{p} :=  \sum_{j\in \mathbb{J}_1(J)}  \| NS(\Lambda_j \phi_j) \|^p_{p}   + \alpha_2 (t, J,n). \label{conver02}
\end{align}
It is easy to see that
\begin{align} \label{conver03}
|\alpha_2 {(t, J,n)}|
& \leq C_{J,p} \sum_{j\neq k, \ j,k \leq \mathbb{J}_1(J)} \int |NS(\Lambda_j \phi_j)(\cdot- x_{j,n})| |NS(\Lambda_k \phi_k)(\cdot- x_{k,n})|^{p-1}     dx.
\end{align}
Define
$$
\displaystyle h_n(t) = \int_{\rt}  |NS(\Lambda_j \phi_j)(\cdot- x_{j,n})| |NS(\Lambda_k \phi_k)(\cdot- x_{k,n})|^{p-1}   dx .
$$
The desired conclusion follows if we obtain that (up to a subsequence for $n\in \mathbb{N}$)
\begin{equation}\label{inequasm}
\lim_{n \to \infty}\sup_{0<t<T}  |h_n(t)| =0 \quad  \textrm{for}\,\, j\neq k.
\end{equation}
Using the fact $NS(\Lambda_j \phi_j)  \in C([0,T], L^p)$,  we can show that $h_n(t)$ is uniformly bounded and
 equi-continuous. In addition, due to the orthogonality of cores, it can be easily checked that  $\lim\limits_{n\to \infty} h_n(t) = 0$, for each
$t \in [0,T]$. Now we apply the Arzel\`{a}-Ascoli Theorem, there exists subsequence $(n_k)_{k \geq 1}$ such that $\displaystyle \lim_{n_k\to \infty} \sup_{0<t<T}h_{n_k}(t) = 0$. Lemma~\ref{Lem:orth} is proved. $\hfill\Box$\\

Now let us turn our attention to the proof of Lemma~\ref{smfor}, which shows the smallness of the forcing term. First we present the main inequalities that will be used, whose proof relies heavily on the space-time estimate of heat kernel.

\begin{prop}\label{impineq}
Let $ 3 < p < \infty$,    $1/\beta(p) = 1/2 -  3/10p$.
\begin{itemize}
\item[\rm (i)] Let $1/m = 1/p -\e, \, 0< \e \leq 2/5p $.  We have
 $$\displaystyle \|v   \etl u_0\|_{L_T^{\be}L_x^{\al}} \leq
C T^{  (1 - 3/m)/2 }\|u_0\|_{L_x^m} \|v\|_{L_T^{5p/3}L_x^{5p/3}}.
$$
\item[\rm (ii)] Let $1/r = 1/p + \varepsilon, \, 0 < \varepsilon < \min\{1/3 - 1/p, 4/5p\}$. We have
$$
\|v   \etl u_0\|_{L_T^{\be}L_x^{\al}} \leq C T^{ (1 - 3/r)/2}\|u_0\|_{L_x^{r}} \|v\|_{L_T^{5p/3}L_x^{5p/3}}.
$$

\item[\rm (iii)] Let $\smp = 4p/3, \, s_{\sigma(p),p} = 3(1/\smp- 1/p)$. We have
$$
\|v  \etl u_0\|_{L_T^{\be}L_x^{\al}} \leq C T^{ (1 - 3/p)/2}\|u_0\|_{\dot B^{s_{\sigma(p),p}}_{\smp,\smp}} \|v\|_{L_T^{5p/3}L_x^{5p/3}}.
$$
\end{itemize}
\end{prop}
{\bf{Proof.}} First, we prove $\mathrm{(i)}$.   Let $\nu(p)$ satisfy $3/m = 9/5p + 2/\nu(p)$. It is easy to see that $p< m \leq \min\{\nu(p), 5p/3\}$.  Applying H\"older's inequality, one has that
$$
\begin{aligned}
\|v \etl u_0\|_{L_T^{\be}L_x^{\al}}
& \leq T^{ (1 - \frac{3}{m})/2 }\|\etl u_0\|_{L_T^{\nu(p)}L_x^{5p/3}} \|v\|_{L_T^{5p/3}L_x^{5p/3}} \\
& \lesssim T^{ (1 - \frac{3}{m})/2}\| u_0\|_{L_x^m} \|v\|_{L_T^{5p/3}L_x^{5p/3}},
\end{aligned}
$$
where we have used \eqref{eq4.1.8aa}.

Next, we show $\mathrm{(ii)}$.  Let $\rho(p)$ satisfy $2/\rho(p)  + 9/5p  = 3/r$. It follows from the hypothesis that
$3 < r < 5p/3, \rho(p)$. Using the same way as in the proof of (i), we have
$$
\begin{aligned}
\|v  \etl u_0\|_{L_T^{\be}L_x^{\al}}
& \leq  T^{ (1 - \frac{3}{r})/2}\|\etl u_0\|_{L_T^{\rho(p)}L_x^{5p/3}} \|v\|_{L_T^{5p/3}L_x^{5p/3}}  \\
& \lesssim T^{ (1 - \frac{3}{r})/2}\|u_0\|_{L_x^{r}} \|v\|_{L_T^{5p/3}L_x^{5p/3}}.
\end{aligned}
$$
Finally, we prove $\mathrm{(iii)}$.  Let $\eta(p) = 3( 1/\smp  - 3/ 5p ), \,  2/\kappa(p)  +
9/5p  = 3/p$. Obviously, $\eta(p) > 0, \kappa(p) > \smp  $. By Proposition~\ref{embedbs}, we see that
$\dot B^{\eta(p)}_{\smp,\smp} \hookrightarrow L^{5p/3}$. So, we have
$$
\begin{aligned}
\|v   \etl u_0\|_{L_T^{\be}L_x^{\al}}
& \lesssim \big\| \|v(t)\|_{L_x^{5p/3}} \|\etl u_0\|_{\dot B^{\eta(p)}_{\smp,\smp}} \big\|_{\ltb} \\
& \leq T^{ (1 - 3/p)/2}\|\etl u_0\|_{L_T^{\kappa(p)}\dot B^{\eta(p)}_{\smp,\smp}} \|v\|_{L_T^{5p/3}L_x^{5p/3}}.
\end{aligned}
$$
 Finally, applying the result of~(\ref{besovheat}), %$\mathrm{(ii)}$ in Proposition~\ref{lineheatesti},
we can bound the above term by
$$
C T^{ (1 - 3/p)/2}\| u_0\|_{\dot B^{\eta(p)- \frac{2}{\kappa(p)}}_{\smp,\smp}} \|v\|_{L_T^{5p/3}L_x^{5p/3}}.
$$
A simple calculation shows $\eta(p)- \frac{2}{\kappa(p)} = 3(\frac{1}{\smp}- \frac{1}{p}) = s_{\sigma(p),p} $. Hence we obtain the desired result. $\hfill\Box$\\

\noindent {\bf{Proof of Lemma~\ref{smfor}:}}  we rewrite $G_n^J$ so as to make use of the smallness of profiles without constant scales.
We make a cut-off on $\phi_j$ and denote
\begin{equation}\label{smpnota}
\begin{aligned}
U^{0}_{n,\eta}: &= \sjjz{0}\lmdst \phi_{j,\eta}, \quad U^{\infty}_{n,\eta}:= \sjjz{\infty}\lmdst \phi_{j,\eta} \\
\psi_{n,\eta}^J :&= \sjjzc\lmdst \phi_{j,\eta^c} + \psi_{n}^J.
\end{aligned}
\end{equation}
Obviously
$$
\U = \sum_{j \in \mathbb{J}_1(J)}NS(\Lambda^{3/p}_{j,n}\phi_j) + \etl(U^{0}_{n,\eta}+U^{\infty}_{n,\eta}+\psi_{n,\eta}^J).
$$
It follows that
$$
G_{n}^{J} = G_{n,\eta}^{J,1} + G_{n,\eta}^{J,2} + G_{n,\eta}^{J,3} + G_{n}^{J,4},
$$
where
$$
\begin{aligned}
G_{n,\eta}^{J,1} &= \big(2\sum_{j \in \mathbb{J}_1(J)}NS(\Lambda^{3/p}_{j,n}\phi_j) + \etl(U^{0}_{n,\eta}+U^{\infty}_{n,\eta}+\psi_{n,\eta}^J)\big)
\otimes \etl U^{0}_{n,\eta} \\
G_{n,\eta}^{J,2} &= \big(2\sum_{j \in \mathbb{J}_1(J)}NS(\Lambda^{3/p}_{j,n}\phi_j) + \etl(U^{0}_{n,\eta}+U^{\infty}_{n,\eta}+\psi_{n,\eta}^J)\big)
\otimes \etl U^{\infty}_{n,\eta}\\
G_{n,\eta}^{J,3} &= \big(2\sum_{j \in \mathbb{J}_1(J)}NS(\Lambda^{3/p}_{j,n}\phi_j) + \etl(U^{0}_{n,\eta}+U^{\infty}_{n,\eta}+\psi_{n,\eta}^J)\big)
\otimes \etl \psi_{n,\eta}^J \\
G_{n}^{J,4} & = \sum_{j,k \in \mathbb{J}_1:j \neq k} \sol{j} \otimes \sol{k}
\end{aligned}
$$
First, we consider the estimate of $G_{n,\eta}^{J,1}$. Let $r $ be as in $\mathrm{(ii)}$ of Proposition~\ref{impineq},  we see that
$$
\begin{aligned}
\|G_{n,\eta}^{J,1}\|_{L_T^{\be}L_x^{\al}}& \leq C T^{ (1 - 3/r)/2}\|U^{0}_{n,\eta}\|_{L_x^{r}} \times \\
&\|2\sum_{j \in \mathbb{J}_1(J)}NS(\Lambda^{3/p}_{j,n}\phi_j) + \etl(U^{0}_{n,\eta}+U^{\infty}_{n,\eta}+\psi_{n,\eta}^J)\|_{L_T^{5p/3}L_x^{5p/3}}.
\end{aligned}
$$
As a result of~(\ref{constbda}) and ~(\ref{linetmbd}),  one has that
$$
\begin{aligned}
\lim_{J \to \infty}\limsup_{n\to \infty} \|2\sum_{j \in \mathbb{J}_1(J)}& NS (\Lambda^{3/p}_{j,n}\phi_j) + \etl(U^{0}_{n,\eta}+U^{\infty}_{n,\eta}+\psi_{n,\eta}^J)\|_{L_T^{5p/3}L_x^{5p/3}} \\
&\lesssim \sum_{j \in \mathbb{J}_1 , \, j \leq J_c} \|NS(\Lambda_{j}\phi_j)\|_{L_T^{\frac{5p}{3}}L_x^{\frac{5p}{3}}} + \sup_{n \geq 1}\|u_{0,n}\|_{p} < \infty.
\end{aligned}
$$
In addition, $r < p$, by $(1)$ of Proposition~\ref{small},
$$
\lim_{n \to \infty} \|U^{0}_{n,\eta}\|_{L_x^{r}} =0
$$
Therefore, we attain
$$
\lim_{J\to \infty}\limsup_{n\to \infty} \|G_{n,\eta}^{J,1}\|_{L_T^{\be}L_x^{\al}} =0.
$$
Next, using  $\mathrm{(i)}$ of Proposition~\ref{impineq}
and $(2)$ of Proposition~\ref{small}, the
estimate of $G_{n,\eta}^{J,2}$ is quite similar to $G_{n,\eta}^{J,1}$ and we omit the details.

Thirdly, we estimate of $G_{n,\eta}^{J,3}$. Let $\smp, s_{\sigma(p),p}$ be as in Proposition~\ref{impineq} $\mathrm{(iii)}$, so the following inequality holds,
$$
\begin{aligned}
\|G_{n,\eta}^{J,3}\|_{L_T^{\be}L_x^{\al}} \leq & \ C T^{\ha (1 - 3/p)}\|\psi_{n,\eta}^J \|_{\dot B^{s_{\sigma(p),p}}_{\smp,\smp}} \\
& \times \|2\sum_{j \in \mathbb{J}_1(J)}NS(\Lambda^{3/p}_{j,n}\phi_j) + \etl(U^{0}_{n,\eta}+U^{\infty}_{n,\eta}+\psi_{n,\eta}^J)\|_{L_T^{5p/3}L_x^{5p/3}}.
\end{aligned}
$$
Fix $q = r = \smp > p, s_{r,p} = s_{\sigma(p),p}$ in Theorem~\ref{thm:lppd}, then it follows from Proposition~\ref{small} $(3)$ that
$$
\lim_{J \to \infty}\lim_{\eta \to \infty}\limsup_{n\to \infty}
\big\|\psi_{n,\eta}^J \big\|_{\dot B^{s_{\sigma(p),p}}_{\smp,\smp}} = 0
$$
so we acquire
$$
\lim_{J \to \infty}\lim_{\eta \to \infty}\limsup_{n\to \infty} \|G_{n,\eta}^{J,3}\|_{L_T^{\be}L_x^{\al}} =0
$$
Finally, the smallness of
$G_{n}^{J,4}$ has been given as in Lemma \ref{Lem:orth}. Combining the estimate on $G_{n,\eta}^{J,1},~G_{n,\eta}^{J,2},~G_{n,\eta}^{J,3}\,~\textrm{and}\,~G_{n}^{J,4}$, we complete the proof of Lemma~\ref{smfor}. $\hfill\Box$

%%%%%%%%%%%%%%%%%%%%%%%%%%%%%%%%%%%%%%%%%%%%%%%%%%%%%%%%%%%%%%%%%%%%%%%%%%%%%%%%%%%%%%%%%%
%%%%%%%%%%%%%%%%%%%%%%%%%%%%%%%%%%%%%%%%%%%%%%%%%%%%%%%%%%%%%%%%%%%%%%%%%%%%%%%%%%%%%%%%%%
%%%%%%%%%%%%%%%%%%%%%%%%%%%%%%%%%%%%%%%%%%%%%%%%%%%%%%%%%%%%%%%%%%%%%%%%%%%%%%%%%%%%%%%%%%

\section{Orthogonality of the profiles of $L^p$-solutions \label{pfthfour}}

This section is intended to prove the orthogonality property of the Navier-Stokes solutions, i.e. formula ~(\ref{soluorth}) in Theorem~\ref{structhm}.
Let us introduce some simplified notations first, define
$$
\begin{aligned}
\anj :&= \sum_{j \in \mathbb{J}_1(J)}NS(\Lambda^{3/p}_{j,n}\phi_j), \ \ \bar{A}^J_n= \sum_{j \in \mathbb{J}_1(J)}NS(\Lambda_j \phi_j) (\cdot-x_{j,n}), \\
\bnj :&=  \lineapp = \etl(U^{0}_{n,\eta}+U^{\infty}_{n,\eta}+\psi_{n,\eta}^J),
\end{aligned}
$$
where $U^{0}_{n,\eta}, U^{\infty}_{n,\eta}, \psi_{n,\eta}^J $ are defined as ~(\ref{smpnota}). Recall that we have
$$
\begin{aligned}
NS(u_{0,n}) = \anj  + \bnj   +R_n^J .
\end{aligned}
$$
{\bf{Proof of ~(\ref{soluorth}).}} Let $0 < t < \widetilde T$ be fixed, using inequality~(\ref{elemineq}) again and the above expression of $NS(u_{0,n})$, we easily find
$$
\begin{aligned}
|\gamma_n^J(t)| &= \big|\|NS(u_{0,n})(t)\|_{ p}^p - \sum_{j \in \mathbb{J}_1(J)} \|NS(\Lambda_{j} \phi_j)(t)\|_{  p}^p - \|\lineapp  \|_{p}^p \big|  \\
& \leq C_p \bigg( \big|\!\int_{\rt}|\anj|^p dx - \sum_{j \in \mathbb{J}_1(J)} \int_{\rt} \big| NS(\Lambda_{j} \phi_j) \big|^p dx\big|  \\
& \qquad \qquad + \int_{\rt} |\rnj|^p dx + \int_{\rt}|\bnj|^{p-1}|\rnj| +|\bnj|\cdot |\rnj|^{p-1} dx    \\
 & \qquad \qquad  + \int_{\rt}|\anj|^{p-1}|\rnj| +|\anj|\cdot|\rnj|^{p-1} dx  \\
 &   \qquad \qquad   +  \int_{\rt}|\anj|^{p-1}|\bnj| +|\anj|\cdot |\bnj|^{p-1} dx  \bigg)
\end{aligned}
$$
it follows from (\ref{infinequorth}) and (\ref{smrem}) that
$$
\begin{aligned}
\lim_{n \to \infty} \big|\!\int_{\rt}|\anj|^p dx &- \sum_{j \in \mathbb{J}_1(J)} \int_{\rt} \big| NS(\Lambda_j \phi_j) \big|^p dx\big| = 0 \\
&\lim_{J\to \infty}\limsup_{n\to \infty}\int_{\rt} |\rnj|^p dx =0
\end{aligned}
$$
Concerning the terms $|\bnj|^{p-1} |\rnj|~ \textrm{and} ~|\bnj|\cdot|\rnj|^{p-1}$, by H\"older's inequality we have
$$
\begin{aligned}
\int_{\rt}(|\bnj|^{p-1}|\rnj| +|\bnj|\cdot |\rnj|^{p-1}) dx  \leq  \|\bnj\|^{p-1}_{p}\|\rnj\|_{p} + \|\bnj\|_{p}\|\rnj\|^{p-1}_{p}.
\end{aligned}
$$
By~(\ref{linetmbd}) and the estimate on $\rnj$, one has that
$$
\lim_{J\to \infty}\limsup_{n\to \infty}  \int_{\rt}(|\bnj|^{p-1}|\rnj| +|\bnj|\cdot |\rnj|^{p-1} )dx =0.
$$
We can deal with the terms $|\anj|^{p-1}|\rnj|~\textrm{and}~|\anj|\cdot|\rnj|^{p-1}$ similarly. Thus we are only left
with terms $|\anj|^{p-1}|\bnj|~\textrm{and}~|\anj|\cdot |\bnj|^{p-1}$. Moreover, in view of Lemma \ref{Lem:orth}, it suffices to show that

\begin{prop}\label{crossest}
Let $\bar{A}_n^J, \bnj $ be as above, $0 < t < \widetilde T$ be fixed, then it holds
$$
\lim_{J\to \infty}\lim_{n\to \infty} \int_{\rt}(|\bar{A}^J_n|^{p-1}|\bnj| +|\bar{A}^J_n|\cdot |\bnj|^{p-1}) dx =0
$$
\end{prop}
{\bf{Proof:}} Let $\epsilon >0$ be arbitrary, as we know
$$
\bnj = \etl(U^{0}_{n,\eta}+U^{\infty}_{n,\eta}+\psi_{n,\eta}^J)
$$
 Concerning $|\bar{A}^J_n|^{p-1}|\bnj|$, it suffices to verify the smallness of $|\bar{A}^J_n|^{p-1}|\etl U^{0}_{n,\eta}|, |\bar{A}^J_n|^{p-1}|\etl U^{\infty}_{n,\eta}| $ and $|\bar{A}^J_n|^{p-1}|\etl \psi_{n,\eta}^J|$. Let us estimate them separately. \\

First, we show that
 $$\displaystyle \lim_{n\to \infty} \int_{\rt} |\bar{A}^J_n|^{p-1}|\etl U^{0}_{n,\eta}| dx = 0.$$
For each $j \in \mathbb{J}_1$, we approximate $NS(\Lambda_j \phi_j)$ by a smooth function with compact support in $\lp$ spaces,
denote this function as $\Theta_{j,t}(x)$, so we see
$$
\begin{aligned}
&\int_{\rt} |\bar{A}^J_n|^{p-1}|\etl U^{0}_{n,\eta}|~ dx \\
& \leq C_{p,J} \sjjz{1} \|NS(\Lambda_j \phi_j)- \Theta_{j,t}\|_{p}^{p-1} \|U^{0}_{n,\eta}\|_{p} +
C_{p,J} (\sjjz{1}\| \Theta_{j,t}\|^{p-1}_{L^{a'(p-1)}})\|U^{0}_{n,\eta}\|_{L^a}
\end{aligned}
$$
where we have used H\"older inequality in the last inequality and $1<a < p, a' ~\textrm{satisfies}~ \frac{1}{a} + \frac{1}{a'} =1$.
 So we have
 $$
 \int_{\rt} |\bar{A}^J_n|^{p-1}|\etl U^{0}_{n,\eta}| dx  \to 0.
 $$
  Similarly, one can show that
$$
\displaystyle \lim_{n\to \infty} \int_{\rt} |\bar{A}^J_n|^{p-1}|\etl U^{\infty}_{n,\eta}| dx = 0.
$$
Finally, we prove that
$$
\displaystyle \lim_{J\to \infty}\lim_{\eta \to \infty}\lim_{n\to \infty} \int_{\rt} |\bar{A}^J_n|^{p-1}|\etl \psi_{n,\eta}^J| dx = 0.
$$
Assume $0< t < T < \widetilde T$, let $ \bar{J} > J_c$  be specified later, $J_c = J_c(T)$ is given by Remark~\ref{rk:5}, choose
$J > \bar{J} $ sufficiently large. Write $\bar{A}^J_n $ into two parts:
$$
\bar{A}^J_n =  \sum_{j \in \mathbb{J}_1 , \, \bar{J} \leq j \leq J }NS(\Lambda_j \phi_j) + \sum_{j \in \mathbb{J}_1 , \, j < \bar{J}}NS(\Lambda_j \phi_j)
$$
By H\"older's inequality,
\begin{equation}\label{resmall}
\begin{aligned}
\int_{\rt} & \big|\sum_{j \in \mathbb{J}_1 , \,  \bar{J} \leq j \leq J  }NS(\Lambda_j \phi_j) \big|^{p-1} |\etl \psi_{n,\eta}^J| dx \\
& \leq \big\| \sum_{j \in \mathbb{J}_1 , \, \bar{J}  \leq j \leq J}NS(\Lambda_j \phi_j)\big\|^{p-1}_{p} \|\psi_{n,\eta}^J\|_{p} \\
&  \lesssim \big( \sum_{j \in \mathbb{J}_1 , \, \bar{J}  \leq j \leq J} \| NS(\Lambda_j \phi_j)\|^p_{p} + \varepsilon(J,n) \big)^{\frac{p-1}{p}} \|\psi_{n,\eta}^J\|_{p}
\end{aligned}
\end{equation}
Here $\lim\limits_{n\to \infty} \varepsilon(J,n) =0$ and we have used ~(\ref{infinequorth}) again. Moreover, for any $J,~\eta$, we have
$$
\begin{aligned}
\limsup_{n\to \infty }\|\psi_{n,\eta}^J\|_{p} &= \limsup_{n\to \infty }\|\sjjzc\lmdst \phi_{j,\eta^c} + \psi_n^J \|_{p}   \\
& \leq \limsup_{n\to \infty } \| \sjjzc\lmdst \phi_{j,\eta^c}\|_{p} + \limsup_{n\to \infty }\|\psi_n^J  \|_{p}.
\end{aligned}
$$
Using ~(\ref{orthod}), we see that  there exists a $N_1=N_1(J,\eta)$, such that for any $n > N_1$, $\displaystyle \|\psi_{n,\eta}^J\|_{p} \leq C\sup_{n \geq 1} \|u_{0,n}\|_{p}$.  Now let $n > N_1$, in view of Remark \ref{rk:5} and noticing that $\bar{J}  > J_c$, one sees that ~(\ref{resmall}) has an upper bound
$$
C  \left((\sum_{j \in \mathbb{J}_1 , \, \bar{J}  \leq j \leq J}  \|  \phi_j \|^{p }_{p} )^{(p-1)/p} + |\varepsilon(J,n)|^{(p-1)/p}\right) \sup_{n \geq 1} \|u_{0,n}\|_{p}
$$
Since  $\sum_{j \geq 1} \|\phi_j\|^p_{p}$ is convergent, we now fix $\bar{J}=\bar{J}(\e) $ such that
\begin{equation}\label{smterone}
C\sum_{j \in \mathbb{J}_1 , \, \bar{J}  \leq j \leq J} \| NS(\phi_j)\|^{p-1}_{p}  \sup_{n \geq 1} \|u_{0,n}\|_{p}  \leq \e.
\end{equation}
Using the fact that $\lim\limits_{n\to \infty} \varepsilon(J,n) =0$, we conclude
there exists a $N_2= N_2(\e,J)$, such that  $n > N_2$
\begin{equation}\label{smtertwo}
C|\varepsilon(J,n)|^{\frac{p-1}{p}} \sup_{n \geq 1} \|u_{0,n}\|_{p} < \e.
\end{equation}
Gathering the estimates~(\ref{smterone}) and~(\ref{smtertwo}), we see that for $n > \max\{N_1,N_2\}$
\begin{equation}\label{ineqtwo}
\int_{\rt} \big|\sum_{j \in \mathbb{J}_1 , \, j \geq \bar{J} }NS(\Lambda_j \phi_j) \big|^{p-1} |\etl \psi_{n,\eta}^J| dx < 2\e.
\end{equation}
Considering the estimate of $\sum\limits_{j \in \mathbb{J}_1 (\bar{J})}NS(\Lambda_j \phi_j)$,   we approximate each $NS(\phi_j)(t,x)$ by $\Theta_{j,t}(x)$
in $\lp$ spaces,
$$
\begin{aligned}
&\int_{\rt} \big|\sum_{j \in \mathbb{J}_1 (\bar{J}) }NS(\Lambda_j \phi_j) \big|^{p-1} |\etl \psi_{n,\eta}^J| dx \\
& \leq C_p \big(\sum_{j\in \mathbb{J}_1 (\bar{J})}\|NS(\Lambda_j \phi_j) - \Theta_{j,t}\|_{p}\big )^{p-1} \|\psi_{n,\eta}^J\|_{p}   \\
& \qquad \qquad \quad +  C_p {t}^{ s_{\smp,p}/2}\| \psi_{n,\eta}^J\|_{\dot B^{s_{\smp,p}}_{\smp,\infty}} \big(\sum_{j\in \mathbb{J}_1 (\bar{J})}
\|\Theta_{t,x}\|_{L^{(p-1)\smp'}} \big)^{p-1}.
\end{aligned}
$$
Here we have used H\"older inequality and equivalent definition of Besov spaces ~(see(\ref{heatbesov})) in deriving the last inequality. From the above,
we know $n > N_1,\, \displaystyle \|\psi_{n,\eta}^J\|_{p} \leq C\sup_{n \geq 1} \|u_{0,n}\|_{p}$, so we choose $\Theta_{j,x}$ sufficiently close to $NS(\phi_j)$ so that
\begin{equation}\label{ineqthree}
 C_p (\sum_{j\in \mathbb{J}_1 (\bar{J})}\|NS(\Lambda_j \phi_j) - \Theta_{j,t}\|_{p}\big )^{p-1} \|\psi_{n,\eta}^J\|_{p} < \e
\end{equation}
while by~(\ref{orthob})~and~Proposition~\ref{small}~$\mathrm{(iii)}$, note we already take $q=r =\smp, s_{r,p} = s_{\smp,p}$ there,
$$
\lim_{J \to \infty}\lim_{\eta \to \infty}\limsup_{n\to \infty}
\big\|\psi_{n,\eta}^J \big\|_{\dot B^{s_{\sigma(p),p}}_{\smp,\smp}} = 0
$$
As $\bar{J}$ is independent of $J,\eta$, $0<t < \widetilde T$ is fixed, so $\exists J(\e), \forall J > J(\e), \,\exists \eta(J,\e), \forall \eta > \eta(J,\e), \,\exists \bar{N}= \bar{N}(\eta,J,\e)$, for $\forall \,n > \bar{N} $
\begin{equation}\label{ineqfour}
 C_p {t}^{\ha s_{\smp,p}}\| \psi_{n,\eta}^J\|_{\dot B^{s_{\smp,p}}_{\smp,\infty}} \big(\sum_{j\in \mathbb{J}_1 (\bar{J})}
\|\Theta_{t,x}\|_{L^{(p-1)\smp'}} \big)^{p-1} < \e
\end{equation}
Gathering~(\ref{ineqtwo}),(\ref{ineqthree})~and~(\ref{ineqfour}), $\exists J(\e), \forall J > J(\e), \,\exists \eta(J,\e), \forall \eta > \eta(J,\e), \,\exists \widetilde N :=\max\{\bar{N},N_1,N_2\}$, s.t for  $\forall\, n > \widetilde N $
$$
\int_{\rt} |\bar{A}^J_n|^{p-1}|\etl \psi_{n,\eta}^J| dx <4\e.
$$
Using the same way as above, one can estimate $|\bar{A}^J_n|\cdot |\bnj|^{p-1}$ and the details are omitted. $\hfill\Box$ \\

%%%%%%%%%%%%%%%%%%%%%%%%%%%%%%%%%%%%%%%%%%%%%%%%%%%%%%%%%%%%%%%%%%%%%%%%%%%%%%%%%%%%%
%%%%%%%%%%%%%%%%%%%%%%%%%%%%%%%%%%%%%%%%%%%%%%%%%%%%%%%%%%%%%%%%%%%%%%%%%%%%%%%%%%%%%
%%%%%%%%%%%%%%%%%%%%%%%%%%%%%%%%%%%%%%%%%%%%%%%%%%%%%%%%%%%%%%%%%%%%%%%%%%%%%%%%%%%%%
\appendix

\section{Perturbation of the Navier-Stokes equation} \label{perturbationns}

We will state some known estimates on Navier-Stokes equations and prove a result on perturbative Navier-Stokes system. First, let us recall two classical inclusions on Besov  and Triebel-Lizorkin spaces, one can refer to Chapter $1$ in \cite{WaHuHaGu11} or Chapter $2$ in~\cite{BaChDa11} for details.
\begin{prop}\label{embedbs}
Assume $ s_1 < s_2 ,\, 1 \leq p_2 < p_1 < \infty,\, 1 \leq q,r \leq \infty, s_1 - \frac{d}{p_1} = s_2 -\frac{d}{p_2}   $, then
$$
   \dot B^{s_2}_{p_2,q} (\R^d)\hookrightarrow \dot B^{s_1}_{p_1,q}(\R^d),  \ \
   \dot F^{s_2}_{p_2,q} (\R^d)\hookrightarrow \dot F^{s_1}_{p_1,r}(\R^d) , \,\, L^{p_1}(\R^d) = \dot F^0_{p_1,2}(\R^d).
$$
\end{prop}
Between Besov and Triebel-Lizorkin spaces, for $s\in \R, 1\leq p < \infty,1 \leq q \leq \infty $, there holds
\begin{equation}\label{rebstr}
\dot B^s_{p,p\wedge q}  \hookrightarrow \dot F^s_{p, q} \hookrightarrow \dot B^s_{p,p\vee q}
\end{equation}
where $ p \wedge q = \min\{p,q\}, \, p\vee q = \max\{p,q \}$.

Now let us recall some known estimates for the NS equation, which were essentially established or hidden in \cite{Ch99,Gi86,Wa04,Pl96} and summarized in \cite{WaHuHaGu11}.

\begin{prop}\label{Prop4.1.2}
Let $a\ge 0$, $1\le r\le p\le \infty$, $0< \lambda\le \infty$ and
$2/\gamma=a+n(1/r-1/p)$.  Then we have
\begin{align}\label{eq4.1.8}
 \|e^{t\Delta} f\|_{L^\gamma(\mathbb{R}_+; \dot B^{0}_{p, \lambda})}\le
C\|f\|_{\dot B^{-a}_{r, \, \lambda \wedge \gamma}}.
\end{align}
In addition, if $\gamma \geq r >1 $, then
\begin{align}\label{eq4.1.8aa}
 \|e^{t\Delta} f\|_{L^\gamma(\mathbb{R}_+; L^p)}\le
C\|f\|_{\dot H^{-a}_r}.
\end{align}
\end{prop}
We remark that \eqref{eq4.1.8aa} is a consequence of \eqref{rebstr}, \eqref{eq4.1.8} and  Proposition \ref{embedbs}.
\begin{prop}\label{Prop4.1.4}
Let $1\le r\le p\le \infty$ and $1< \gamma,
\gamma_1<\infty$ satisfy
\begin{align}\label{eq4.1.22}
\frac{1}{\gamma}= \frac{1}{\gamma_1}+ \frac{k}{2}+
\frac{n}{2}\left(\frac{1}{r}- \frac{1}{p} \right)-1, \ \
\frac{k}{2}+ \frac{n}{2}\left(\frac{1}{r}- \frac{1}{p}\right)<1, \ \
k=0,1.
\end{align}
Then we have
\begin{align}\label{eq4.1.23}
\|\nabla^k\mathscr{A}_{t_0} f\|_{L^\gamma([t_0,\infty); L^p)} \lesssim
\|f\|_{L^{\gamma_1}([t_0,\infty); L^r)}.
\end{align}
\end{prop}

\begin{prop}\label{Prop4.1.6}
Let $1\le r \le \infty$, $1\le q'\le \lambda \le \infty$\footnote{$p'$ stands for the conjugate number of $p$, i.e. $1/p+1/p'=1$.}. Then
\begin{align}
\|\mathscr{A}_{t_0} f\|_{L^\infty([t_0,\infty),  \dot B^{0}_{r, \lambda})}
\lesssim \|f\|_{L^{q'}([t_0,\infty), \dot B^{-2/q}_{r, \lambda})}.
\label{eq4.1.26}
\end{align}
\end{prop}
As a direct consequence, we have
 we now give several estimates adapted to our needs.

\begin{cor}\label{lineheatesti}
Let $s\in \R,3 < p< \infty, \,1 \leq q \leq \gamma \leq \infty,\, 1 \leq m \leq \infty$, $0< T < \infty, \, \frac{1}{\beta(p)} = \frac{1}{2} - \frac{3}{10p}$, then we have
\begin{align}
 & \|\etl u_0\|_{\X} \lesssim  \|u_0\|_{p},  \\
 & \|\etl u_0\|_{L^\gamma ((0,\infty), \dot B^{s+ \frac{2}{\gamma }}_{m,q})} \lesssim \|u_0\|_{\dot B^{s}_{m,q}}. \label{besovheat} \\
 & \| \mathscr{A}_{0}\mathbb{P} \nabla  f \|_{\X} \lesssim  \| f\|_{L_T^{\be}L_x^{\al}} \lesssim  T^{  (1 - \fp)/2}\| f\|_{L_T^{ 5 p/6 }L_x^{5p/6}},
\end{align}
\end{cor}
 {\bf Proof.} By Proposition  \ref{Prop4.1.4} and H\"older's inequality,
$$
\| \mathscr{A}_{0} \mathbb{P} \nabla f \|_{L_T^{5p/3}L_x^{5p/3}}  \lesssim \| f\|_{L_T^{\be}L_x^{\al}}\lesssim  T^{  (1 - \fp)/2}\| f\|_{L_T^{ 5 p/6 }L_x^{5p/6}} .
$$
By Propositions  \ref{Prop4.1.6} and  \ref{embedbs},  taking $s(p) = 1 - 2/\be  >0 $,    we see that
\begin{equation}\label{estispacetimeb}
\begin{aligned}
\| \mathscr{A}_{0}\mathbb{P} \nabla  f \|_{L_T^{\infty}L_x^{p}} \lesssim \| \mathscr{A}_{0}\mathbb{P} \nabla f \|_{L_T^{\infty} \dot B^{s(p)}_{\al,\al}}
 \leq  \|f\|_{L_T^{\be} \dot B^{0}_{\al,\al}} \lesssim \|f\|_{L_T^{\be}L_x^{\al}}.
\end{aligned}
\end{equation}

Next, we perform the perturbation analysis for the Navier-Stokes equation.  Our approach follows  Theorem $3.1$ in~\cite{GaIfPl03}. Similar conclusions in critical spaces can also be found in~\cite{GaKoPl13}.

\begin{prop}\label{perturprop}
Let $ 3<p<\infty,\,w\in \X$ be vector fields in $\rt$, $f \in L_T^{\be}L_x^{\al}$ is a $3\times 3$ matrix function. $\be$ is the same as $\mathrm{Lemma~\ref{bddri}}$, $0< T<\infty$,  $v$ satisfies the following perturbative NS equation (PNS)
$$
 \quad \left\{
\begin{array}{ll}
\partial_t v & + (v\cdot \nabla) v -\Delta v +(v\cdot \nabla) w+ (w\cdot \nabla) v=  -\nabla p  -\nabla\cdot f  \\
\mathrm{div} \, v&=  0 \\
v_{|t=0} &=  v_{0}\, .
\end{array}
\right.
$$
Assume that there exist two constants $\e_0 \ll 1$ and $C \gg 1$ such that
$$
\|v_0\|_{p} + \|f\|_{L_T^{\be}L_x^{\al}} \leq \e_0\,T^{- (1-\fp)/2}\, \mathrm{exp}(-C \,T^{ 5p(1-\fp)/6}\|w\|^{5p/3}_{L_T^{5p/3}L_x^{5p/3}}).
$$
Then $v$ belongs to $\X$, moreover we have
$$
\|v\|_{\X} \lesssim (\|v_0\|_{p} + \|f\|_{L_T^{\be}L_x^{\al}} )\, \exp(C\,T^{ 5p(1-\fp)/6}\|w\|^{5p/3}_{L_T^{5p/3}L_x^{5p/3}})
$$
\end{prop}

{\bf{Proof:}} The Duhamel formula of  (PNS)  is
$$
v = \etl v_0 - \mathscr{B}(v,v) - \mathscr{B}(v,w) - \mathscr{B}(w,v) - \mathscr{A}_0\mathbb{P}\nabla \cdot f
$$
For any time interval $(\alpha,\beta) \subset (0,T)$,  we denote  $\mathfrak{X}^p_{(\alpha,\beta)}: = L^{5p/3}((\alpha,\beta),L_x^{5p/3} )\cap C([\alpha,\beta], L^p_x) $. It follows from Corollary \ref{lineheatesti} that
\begin{equation}\label{keyesti}
\begin{aligned}
\|v\|_{\mathfrak{X}^p_{(\alpha,\beta)}}&\leq L (\|v(\alpha)\|_{p} + \|f\|_{L_T^{\be}L_x^{\al}} )     %{L^{\infty}((\alpha,\beta),L^p)\cap L^{5p/3}((\alpha,\beta),L^{5p/3}) }
 + L \,T^{ (1-\fp)/2} \Big(\|v\|^2_{L^{5p/3}((\alpha,\beta), L_x^{5p/3} )}   \\
&\qquad \qquad \qquad\qquad\qquad+ \|v\|_{L^{5p/3}((\alpha,\beta),L_x^{5p/3})} \|w\|_{L^{5p/3}((\alpha,\beta),L_x^{5p/3})} \Big)
\end{aligned}
\end{equation}
for some constant $ L> 1$. Based on this estimate, we can solve (PNS) locally and then extend  to the maximal solution, i.e. solution with maximal existence time $ T (v_0)$.  By the absolute continuity of Lebesgue integral, there exists $N+1$ real numbers $ (T_i)_{0\leq i\leq N} $ such that
$T_0 =0$ and $T_N= T $, satifying $[0,T]= \cup_{i=0}^{N-1}[T_i,T_{i+1}]$ and
\begin{equation}\label{inequalitya}
\begin{aligned}
 \frac{1}{8L} \leq T^{ (1-\fp)/2}&\|w\|_{L^{5p/3}((T_i,T_{i+1}),L_x^{5p/3})} \leq \frac{1}{4L} \quad \forall i \in\{0,\ldots,N-2\} \\
T^{ (1-\fp)/2}&\|w\|_{L^{5p/3}((T_{N-1},T),L_x^{5p/3})} \leq \frac{1}{4L}.
\end{aligned}
\end{equation}
By virtue of~(\ref{inequalitya}), one sees that $N$ satisfying
\begin{equation}\label{estin}
\|w\|^{5p/3}_{L_T^{5p/3}L_x^{5p/3}} \geq \sum_{i=0}^{N-2}\int_{T_i}^{T_{i+1}} \|w\|^{5p/3}_{L_x^{5p/3}} dt \geq
(N-1) \big(\frac{1}{8LT^{ (1-\fp)/2}}\big)^{5p/3}.
\end{equation}
Now let us assume that $\|v_0\|_{p} + \|f\|_{L_T^{5p/3}L_x^{5p/3}}$ is bounded by
\begin{equation}\label{ineqinidata}
 (8L T^{  (1-\fp)/2} (2L)^5)^{-1}
\exp \big[- (8L T^{ (1-\fp)/2})^{5p/3}\|w\|^{5p/3}_{L_T^{5p/3}L_x^{5p/3}} \ln(2L) \big ].
\end{equation}
By continuity in time, we define a maximal time $ \bar{T} < T (v_0)$ such that
\begin{equation}\label{inequalityb}
T^{ (1-\fp)/2} \|v\|_{L^{5p/3}((0,\bar{T}),L^{5p/3})} \leq \frac{1}{4L}.
\end{equation}
If $\bar{T} \geq T$, then by~(\ref{keyesti}), we can   get desired conclusion. \\
\indent If $\bar{T} < T$, then there exists $k \in \{0,\ldots, N-1\}$, such that $T_k \leq \bar{T} < T_{k+1}$. Inserting ~(\ref{inequalitya}) and~(\ref{inequalityb}) into~(\ref{keyesti}), we can get for any $1 \leq i < k$
$$
\|v\|_{\mathfrak{X}^p_{(T_i,T_{i+1})}} \leq L(\|v(T_i)\|_{p} + \|f\|_{L_T^{5p/3}L_x^{5p/3}}) +
2^{-1} \|v\|_{L^{5p/3}((T_i,T_{i+1}),L_x^{5p/3})}.
$$
Since $\mathfrak{X}^p_{(T_i,T_{i+1})}= L^{5p/3}((T_i,T_{i+1}),L_x^{5p/3} )\cap C([T_i,T_{i+1}], L^p_x) $, so we have

\begin{equation}\label{inequalityc}
\|v(T_{i+1})\|_{p} \leq \|v\|_{\mathfrak{X}^p_{(T_i,T_{i+1})}} \leq 2L(\|v(T_i)\|_{p} + \|f\|_{L_T^{5p/3}L_x^{5p/3}}).
\end{equation}
Iterating the above inequality~(\ref{inequalityc}) with respect to $i,\, 1\leq i < k$
\begin{equation}\label{inequalityd}
\|v(T_i)\|_{p} \leq (2L)^{i+1} (\|v_0\|_{p} + \|f\|_{L_T^{5p/3}L_x^{5p/3}}).
\end{equation}
By~(\ref{inequalityd}) and~(\ref{inequalityc}), we see
$$
\|v\|_{\mathfrak{X}^p_{(T_i,T_{i+1})}} \leq 2(2L)^{i+2}(\|v_0\|_{p} + \|f\|_{L_T^{5p/3}L_x^{5p/3}}).
$$
Similarly we have
$$
\|v\|_{\mathfrak{X}^p_{(T_k,\bar{T})}} \leq 2(2L)^{k+2}(\|v_0\|_{p} + \|f\|_{L_T^{5p/3}L_x^{5p/3}}).
$$
Further,
$$
\begin{aligned}
\|v\|_{L_{\bar{T}}^{5p/3}L_x^{5p/3}}
& \leq \sum_{i=0}^{k-1} \|v\|_{L^{5p/3} ((T_i,T_{i+1}),L_x^{5p/3}) }+ \|v\|_{L^{5p/3}((T_k, \bar{T}),L_x^{5p/3})} \\
& \leq (2L)^{N+4}(\|v_0\|_{p} + \|f\|_{L_T^{5p/3}L_x^{5p/3}}).
\end{aligned}
$$
By~(\ref{estin}) and~(\ref{ineqinidata}), we  have
$$
\|v_0\|_{p} + \|f\|_{L_T^{5p/3}L_x^{5p/3}} \leq \frac{1}{8L(2L)^{N+4}  T^{\ha (1-\fp)}}
$$
with this we finally obtain
$$
T^{\ha (1-\fp)} \|v\|_{L_{\bar{T}}^{5p/3}L_x^{5p/3}} \leq  1/8L,
$$
which contradicts the maximality of $\bar{T}$. Therefore, $\bar{T} \geq T$, and using the above iteration, we have
$$
\|v\|_{\X} \leq (2L)^{N+4}(\|v_0\|_{p} + \|f\|_{L_T^{5p/3}L_x^{5p/3}}). \\
$$
Together with the estimate~(\ref{estin}), we can rewrite the above estimate as
$$
\|v\|_{\X} \leq C  (\|v_0\|_{p} + \|f\|_{L_T^{5p/3}L_x^{5p/3}}) \exp(C T^{   5p(1-\fp)/6} \|w\|^{5p/3}_{L_T^{5p/3}L_x^{5p/3}} )
$$
for some constant  $C\gg 1$. Hence we complete the proof. $\hfill\Box$\\

As a direct consequence, we have the following corollary.
\begin{cor}\label{lscontinuous}
The mapping $w_0  \mapsto T(w_0)$ is lower semi-continuous for (NS) with initial in $L^p,\, 3< p < \infty$, i.e. for arbitrary $\varepsilon >0$, there
exists a $\eta >0$, such that if $\|v_0\|_{p} < \eta$, then $ T(w_0 + v_0) \geq T(w_0)- \varepsilon $.
\end{cor}

\noindent {\bf Sketch Proof of Remark \ref{rk:5}.} Considering the map
\begin{align} \label{mappingapp}
\mathscr{T}: \ u(t) =  e^{ t \Delta } u_0 +  \mathscr{B} (u, u),
\end{align}
by Corollary \ref{lineheatesti}, we have
\begin{align} \label{mappingapp}
\| \mathscr{T}  u \|_{\mathfrak{X}^p_T}  & \leq    \|e^{ t \Delta } u_0\|_{\mathfrak{X}^p_T} + \| \mathscr{B} (u, u)\|_{\mathfrak{X}^p_T} \nonumber\\
&  \lesssim     \|  u_0\|_{p} + T^{(1-3/p)/2}\|u\otimes u\|_{L^{5p/6}_TL^{5p/6}_x} \nonumber\\
& \leq    \|  u_0\|_{p} + T^{(1-3/p)/2} \|u \|^2_{\mathfrak{X}^p_T }.
\end{align}
So, if
$$
u\in \mathscr{D}_T := \{u\in \mathfrak{X}^p_T: \ \|u\|_{\mathfrak{X}^p_T} \leq 2C \|  u_0\|_{p} \},
$$
we can choose $T$ satisfying $ T^{(1-3/p)/2} \|  u_0\|_{p} \leq 1 /4C^2$. It follows that $\mathscr{T}  u \in \mathscr{D}_T$. Using contraction mapping principle, we obtain that NS has a unique solution $u\in \mathfrak{X}^p_T$ and the above proof implies the result of Remark \ref{rk:5}. $\hfill \Box$ \\

\noindent{\bf Acknowledgment.} The work was
supported in part by the National Science Foundation of China, grants   11271023. The second and third named authors were supported in part by Mathematics and Physics Unit ``Multiscale Analysis Moddelling and Simulation'', Top Global University Project, Waseda University. Part of the paper was carried out when the first named author was visiting Laboratoire J.~A.~Dieudonn\'e and he is grateful to Professor Fabrice Planchon for his valuable suggestions and comments and also the support from the China
Scholarship Council.


\medskip
\footnotesize
\begin{thebibliography}{100}

\bibitem{Ab15} B. Abe,  The Navier¨CStokes equations in a space of bounded
functions, Commun. Math. Phys. \textbf{338}(2015), 849--865.

\bibitem{AuDuTc06} P. Auscher, S. Dubois and P. Tchamitchian,  On the stability of global solutions to Navier--Stokes equations in the space, J. Math. Pures Appl., \textbf{83} (2004) 673¨C697

\bibitem{BaBiTa12} H. Bae, A. Biswas and E. Tadmor,  Analyticity and decay estimates
of the Navier--Stokes equations in critical
Besov spaces, Arch. Rational Mech. Anal., {\bf 205} (2012), 963--991.



\bibitem{BL} J. Bergh and J. L\"{o}fstr\"{o}m,
 Interpolation Spaces,   Springer--Verlag,  1976.

\bibitem{BoPa08} J. Bourgain and N. Pavlovic,  Ill-posedness of the Navier--Stokes equations in a critical
space in 3D. J. Funct. Anal., $\mathbf{255}$ (2008), 2233--2247.

\bibitem{BaChDa11}
H.~Bahouri, J.~Chemin and R.~Danchin.
\newblock {\em Fourier analysis and nonlinear partial differential equations}.
\newblock Grundlehren der mathematischen Wissenschaften, 343. Springer, Heidelberg, 2011.

\bibitem{Ca95} M. Cannone,  Ondelettes, paraproduits et Navier-Stokes. Diderot Editeur, Paris (1995).

\bibitem{CaMe95}  M. Cannone, Y. Meyer,  Littlewood-Paley decomposition and Navier-Stokes equations, Methods Appl.
Anal. \textbf{2} (1995), 307--319.

\bibitem{Can97} M. Cannone, A generalization of a theorem by Kato on Navier--Stokes equations, Rev. Mat. Iberoamericana, 13 (1997) 515--541.


\bibitem{Ch99} J.-Y. Chemin,  Th\'{e}or\`{e}mes d'unicit\'{e} pour le syst\`{e}me de Navier--Stokes tridimensionnel,
J. Anal. Math., {\bf 77} (1999), 27--50.

\bibitem{ChGaPa11} J.Y. Chemin, I. Gallagher, and M. Paicu, Global regularity for some classes of large solutions to the Navier--Stokes equations, Ann. of Math.,  \textbf{173} (2011), 983--1012.

\bibitem{CoMoPi14} J. C. Cortissoz, J. A. Montero and C. E. Pinilla,
On lower bounds for possible blow-up solutions to the periodic Navier-Stokes equation
J. Math. Phys. \textbf{55} (2014), 033101.

\bibitem{DoDu09} H. Dong and D. Du, The Navier-Stokes equation in the critical Lebesgue space, Commun. Math. Phys., {\bf 292} (2009), 811-827.


\bibitem{EsSeSv03} L. Escauriaza, G. Seregin and V. Sverak, $L_{3,\infty}$ solutions of Navier--Stokes equations
and backward uniquness, Uspekhi Mat. Nauk., \textbf{58} (2003), 3--44.



\bibitem{FoTe89}  C. Foias and R. Temam,  Gevrey class regularity for the solutions of the Navier--Stokes
equations. J. Funct. Anal., \textbf{87} (1989), 359--369.

\bibitem{GaIfPl03}
I.~Gallagher, D.~Iftimie and F.~Planchon.
\newblock Asympototics and stability for global solutions to the Navier-Stokes equations.
\newblock{\em Ann.Inst.Fourier(Grenoble)}, 53(5):1387--1424, 2003.


\bibitem{GaKoPl13} I. Gallagher, G. S. Koch and F. Planchnon, A profile decomposition approach to the $L^\infty_t(L^3_x)$ Navier-Stokes regularity criterion, Math. Ann. {\bf 355} (2013), 1527--1559.

\bibitem{Ge08} P. Germain,  The second iterate for the Navier--Stokes equation. J. Funct. Anal., $\mathbf{255}$ (2008), 2248--2264.



\bibitem{Gi86} Y. Giga,   Solutions for semilinear parabolic equations in $L_p$ and regularity of weak solutions of the Navier-Stokes system,
J. Differ. Equations \textbf{62} (1986), 182--212.

\bibitem{GiInMa99} Y. Giga, K. Inui and S.  Matsui, On the Cauchy problem for the Navier-Stokes equations with nondecaying
initial data. In: Advances in fluid dynamics, vol. 4 of Quad.Mat., pp. 27-68. Dept.Math., Seconda Univ.
Napoli, Caserta (1999).

\bibitem{GiMi85} Y. Giga and T. Miyakawa, Solutions in $L^r$   of the Navier-Stokes initial value problem,  Arch. Rational Mech. Anal.  {\bf 89}  (1985),  267--281.

\bibitem{GiMi89} Y. Giga and T. Miyakawa, Navier-Stokes flow in $\mathbb{R}^3$  with measures as initial vorticity and Morrey spaces, Comm. Partial Differential Equations,   {\bf 14}  (1989),  577--618.

\bibitem{Ge98} P. G\'{e}rard, Description du d\'{e}faut de compacit\'{e} de l'injection de Sobolev, ESAIM Control Optim. Calc. Var. \textbf{3} (1998), 213--233.

\bibitem{HuWa13} C.  Huang  and  B.  Wang,  Analyticity for the (generalized) Navier-Stokes equations with rough initial data, arXiv:1310.2141.

\bibitem{Iw10} T. Iwabuchi, Navier--Stokes equations and nonlinear heat equations
in modulation spaces with negative derivative indices, J. Differential Equations, {\bf 248} (2010), 1972--2002.

\bibitem{Ka84} T. Kato, Strong $L^p$ solutions of the Navier--Stokes equations in $ \mathbb{{R}}^m$ with applications
to weak solutions,  Math. Z., \textbf{187} (1984), 471--480.

\bibitem{KeKo11} C. E. Kenig, G. S. Koch,  An alternative approach to regularity for the Navier--Stokes equations in critical spaces, Ann.  l'Inst. H. Poincare (C) Non Linear Anal.,   \textbf{28} (2011),   159-187.

\bibitem{KeMe08} C. E. Kenig, F. Merle,  Global well-posedness, scattering and blow-up for the energy critical focusing nonlinear wave equations, Acta Math.,   \textbf{201} (2008),   147--212.


\bibitem{KoTa01} H. Koch, D. Tataru, Well-posedness for the Navier--Stokes equations, Adv. Math., {\bf 157} (2001) 22--35.

\bibitem{Ko10}
G. S.~Koch,
\newblock Profile decompositions for critical Lebesgue and Besov space embeddings.
\newblock{\em Indiana Univ. Math.J.}, \textbf{59} (2010), 1801--1830.

\bibitem{KoNaSeSv09} G. S.~Koch, N. Nadirashvili, G.A. Seregin, V. Sverak, Liouville theorems for the Navier-Stokes
equations and applications, Acta Math., \textbf{203} (2009), 83-105.

\bibitem{KoOg02} H. Kozono, T. Ogawa, Y. Taniuchi, The critical Sobolev inequalities in Besov
spaces and regularity criterion to some semi-linear evolution equations, Math.
Z., \textbf{242} (2002), 251--278.


\bibitem{Le34} {J. Leray,} Sue le mouvemeblt d'ubl liouide visoueux emplissant l'espace, Acta Math., {\bf 63} (1934), 193--248.


\bibitem{Pl96} F. Planchon, Global strong solutions in Sobolev or Lebesgue spaces to the incompressible Navier--Stokes equations in $\mathbb{R}^3$,
Ann. Inst. H. Poincare, AN, {\bf 13} (1996), 319--336.

\bibitem{planchon}
F. Planchon.
\newblock Asymptotic behavior of global solutions to the {N}avier - {S}tokes
equations in {$\rt$}.
\newblock {\em Rev. Mat. Iberoamericana}, 14(1):71--93, 1998.


\bibitem{PoRaSiTi94} G. Ponce, R. Racke, T. C. Sideris, and E. S. Titi, Global stability of large
solutions to the 3D Navier-Stokes equations, Comm. Math. Phys. \textbf{159} (1994),
329--341.

\bibitem{Po15}
E.~Poulon.
\newblock About the possibility of minimal blow up for Navier-Stokes
  solutions with data in {$\dot{H}^s(\R^3)$}.
\newblock{\em arXiv : 1505.06197}, May 2015.

\bibitem{Se12} G. Seregin, A certain necessary condition of potential
blow up for Navier-Stokes equations, Comm. Math. Phys., {\bf 312} (2012), 833--845.

\bibitem{RoSaSi12} J. C. Robinson, W. Sadowski  and R. P. Silva,  Lower bounds on blow up solutions of the three dimensional Navier-Stokes equations
in homogeneous Sobolev spaces, Journal of Mathematical Physics, {\bf 53}, (2012) 115618-1--15.


\bibitem{Tr83} H. Triebel,  {\it Theory of Function Spaces,}
  Birkh\"{a}user--Verlag,  1983.

\bibitem{Wa04} B. Wang, Exponential Besov spaces and their applications to certain evolution equations with
dissipations, Commun. Pure Appl. Anal., \textbf{3} (2004), 883--919.

\bibitem{WaHuHaGu11} B.  Wang, Z.   Huo, C.  Hao and Z.  Guo, \newblock {\em Harmonic analysis method for nonlinear evolution equations. I}.
\newblock{ World Scientific Publishing Co. Pte. Ltd.}, Hackensack, NJ, 2011.

\bibitem{Wa06}  B.  Wang, L.  Zhao, B. Guo,   {\rm  Isometric decomposition operators, function spaces $E^\lambda_{p,q}$ and their applications to nonlinear evolution equations,}  J. Funct. Anal., {\bf 233} (2006), 1--39.

\bibitem{Wa15} B.  Wang, Ill-posedness for the Navier-Stokes equation in critical Besov spaces $\dot B^{-1}_{\infty, q}$, Adv. in Math., \textbf{268} (2015), 350--372.

\bibitem{We80}    F. B. Weissler, The Navier-Stokes initial value problem in Lp, Arch. Rational
Mech. Anal. \textbf{74} (1980), 219--230.


\bibitem{Yo10} T. Yoneda,  Ill-posedness of the 3D Navier--Stokes equations in a generalized Besov
space near $BMO^{-1}$. J. Funct. Anal., $\mathbf{258}$ (2010), 3376--3387.



























\end{thebibliography}
\end{document}